\documentclass[english]{article}
\usepackage[T1]{fontenc}
\usepackage[letterpaper]{geometry}
\usepackage{float}
\usepackage{amsmath}
\usepackage{graphicx}
\usepackage{color}
\usepackage{epstopdf}

\usepackage[colorlinks, hypertexnames=false]{hyperref}
\hypersetup{linkcolor=blue, citecolor=red} % set color for references

\usepackage{adjustbox}
\usepackage{bm}
\usepackage{multirow}
\usepackage{array}
\usepackage{mathrsfs}
\usepackage{amsmath, amssymb, amsthm}
\newtheorem{proposition}{Proposition}
\usepackage{subfigure}
\usepackage{tikz}
\usepackage{algorithm}
\usepackage{algpseudocode}
\newtheorem{definition}{Definition}

\newtheorem{remark}{Remark}

\usepackage{amsfonts} 
\usepackage{color}
\usepackage{adjustbox}
\definecolor{black}{rgb}{0,0,0}

\definecolor{red}{rgb}{1,0,0}

\definecolor{blue}{rgb}{0,0,1}

\usepackage{multirow}

\makeatother
\usepackage{mathtools}
\usepackage{esint}
\usepackage{cleveref}
\usepackage{babel}
\usepackage{authblk}
\newcommand{\bfk}{{\bf k}}
\newcommand{\bfx}{{\bf x}}
\def\d{\mathrm{d}}

\title{}

\title{\textbf{}}
\title{\textbf{Uncertainty Quantification of Nonlinear Lagrangian Data Assimilation Using Linear Stochastic Forecast Models}}
\author{Nan Chen\thanks{Department of Mathematics, University of Wisconsin-Madison, WI, USA.  (chennan@math.wisc.edu) }, \quad
	Shubin Fu\thanks{Eastern Institute for Advanced Study, Eastern Institute of Technology, Ningbo, Zhejiang 315200, P. R. China.  Corresponding author (shubinfu89@gmail.com, shubinfu@eias.ac.cn).} \;
	
}

\begin{document}
	\maketitle
	\begin{abstract}
Lagrangian data assimilation exploits the trajectories of moving tracers as observations to recover the underlying flow field. One major challenge in Lagrangian data assimilation is the intrinsic nonlinearity that impedes using exact Bayesian formulae for the state estimation of high-dimensional systems. In this paper, an analytically tractable mathematical framework for continuous-in-time Lagrangian data assimilation is developed. It preserves the nonlinearity in the observational processes while approximating the forecast model of the underlying flow field using linear stochastic models (LSMs). A critical feature of the framework is that closed analytic formulae are available for solving the posterior distribution, which facilitates mathematical analysis and numerical simulations. First, an efficient iterative algorithm is developed in light of the analytically tractable statistics. It accurately estimates the parameters in the LSMs using only a small number of the observed tracer trajectories. Next, the framework facilitates the development of several computationally efficient approximate filters and the quantification of the associated uncertainties. A cheap approximate filter with a diagonal posterior covariance derived from the asymptotic analysis of the posterior estimate is shown to be skillful in recovering incompressible flows. It is also demonstrated that randomly selecting a small number of tracers at each time step as observations can reduce the computational cost while retaining the data assimilation accuracy. Finally, based on a prototype model in geophysics, the framework with LSMs is shown to be skillful in filtering nonlinear turbulent flow fields with strong non-Gaussian features.
	
\end{abstract}


\section{Introduction}
Lagrangian data assimilation exploits the trajectories of moving tracers, such as drifters or floaters, as observations to recover the underlying flow field \cite{apte2013impact, apte2008data, apte2008bayesian, ide2002lagrangian}. It is essential for state estimation and prediction in geophysics, climate science, and hydrology \cite{griffa2007lagrangian, blunden2019look, honnorat2009lagrangian, salman2008using, castellari2001prediction}, especially in situations that the traditional Eulerian observations are not available. Critical applications include the Global Drifter Program \cite{centurioni2017global}, which aims to estimate near-surface currents by tracking the surface drifters deployed throughout the global ocean, and the Argo program \cite{gould2004argo} that utilizes a fleet of robotic instruments drifting with the ocean currents for advancing the operational ocean data assimilation. Lagrangian data assimilation also provides a powerful tool to facilitate the recovery of the ocean eddies in the Arctic regions, where the sea ice floes play the role of the Lagrangian tracers \cite{mu2018arctic, chen2022efficient}.

Despite its unique advantages in facilitating the state estimation of many practical problems, there are several challenges in applying Lagrangian data assimilation to complex systems.
First, one of the fundamental difficulties in Lagrangian data assimilation is the intrinsic nonlinearity \cite{apte2008bayesian, chen2014information}, which comes from not only the governing equations of the underlying flow field but the observational processes as well. This is quite different from the traditional data assimilation with Eulerian observations, in which the observations are often given by a linear combination of the state variables. The observational processes in Lagrangian data assimilation appear as a set of dynamical equations based on Newton's law that relates the time derivative of the tracer displacement with the flow velocity. Since the velocity field is usually a complicated nonlinear function of the tracer locations, the overall observational processes are highly nonlinear. As a result, adopting the exact Bayesian formula for state estimation becomes almost impractical. The intrinsic nonlinearity in the observational processes also requires a careful design of suitable numerical schemes for data assimilation to reduce the filtering error and prevent filter divergence \cite{apte2013impact, sun2019lagrangian, salman2008hybrid}.
Second, the underlying flow field that drives the Lagrangian tracers is typically high-dimensional with multiscale features \cite{pedlosky1987geophysical, vallis2017atmospheric}. A high-resolution numerical solver is needed to integrate the complex forecast model to resolve the dynamics across different spatiotemporal scales and guarantee numerical stability. However, such a procedure is often computationally expensive, especially when it is incorporated into the ensemble data assimilation framework that requires repeating the costly numerical integration multiple times in each assimilation cycle.
Third, the forecast model of the underlying flow field in Lagrangian data assimilation is typically built under the Eulerian coordinate. Therefore, the computational burden is further increased as a result of the coordinate transformation that has to be carried out constantly in the entire data assimilation procedure. Interpolating the Lagrangian data to fixed Eulerian mesh grids also introduces extra errors and affects the accuracy of the recovered states \cite{mead2005assimilation, honnorat2010identification}.

To reduce the computational cost, suitable surrogate forecast models are often utilized to replace the complicated physical system in data assimilation, aiming to provide an accurate short-term forecast of the probability density function (PDF) \cite{law2015data, majda2018model, gershgorin2010improving, delsole2004stochastic}. Notably, these approximate forecast models do not necessarily need to capture the exact physics of the original system, as their primary role is the statistical forecast. Thus, a significant degree of freedom is allowed in developing these surrogate models. Among different approximate models, linear stochastic models (LSMs) are widely adopted in solving complex data assimilation problems, where stochastic noise is utilized to effectively parameterize the nonlinear deterministic time evolution of many chaotic or turbulent dynamics \cite{majda2016introduction, farrell1993stochastic, berner2017stochastic, branicki2018accuracy, majda2018model, li2020predictability, harlim2008filtering, kang2012filtering}.
While data assimilation with LSMs may lead to significant errors in recovering strong non-Gaussian features, including extreme events, in the situation with sparse Eulerian observations \cite{majda2012filtering, majda2010mathematical}, the unique properties of the multiple moving observations can facilitate the use of LSMs in recovering the underlying flow field in Lagrangian data assimilation. The overall complexity of  Lagrangian data assimilation is significantly reduced by virtue of LSMs, as the remaining nonlinearity comes merely from the observational processes. Notably, although the observational processes involve strong nonlinear interactions between the flow and the tracer variables, these processes are often linear functions of the flow variables conditioned on the given tracer displacements \cite{chen2014information, apte2013impact, slivinski2015hybrid}. Such conditional linearity plays a crucial role in developing closed analytic formulae for solving the posterior distribution in Lagrangian data assimilation \cite{chen2014information, chen2015noisy}. On the other hand, unlike data assimilation with sparse Eulerian observations, the multiple Lagrangian tracers can effectively make the posterior distribution trust more towards the information provided by the observations. This advances the recovery of extreme events and intermittency even with LSMs. In addition, the LSMs of the underlying flow field can often be written into the form of its spectral representation, which avoids the constant and expensive transformation between Eulerian and Lagrangian coordinates and automatically eliminates the interpolation errors. Given the potential prospects of utilizing LSMs to facilitate Lagrangian data assimilation, it is of practical importance to understand the associated uncertainties in such a strategy.

In this paper, an analytically tractable mathematical framework for continuous-in-time nonlinear Lagrangian data assimilation is developed. It advances the study of the uncertainty quantification (UQ) of nonlinear Lagrangian data assimilation with LSMs.
First, parameter estimation is the prerequisite for data assimilation. Since the underlying flow field is not observed directly, estimating the parameters in the LSMs for describing the flow field has to be carried out using the observed tracer trajectories. Typical parameter estimation methods {include  traditional approaches}  such as the Markov chain Monte Carlo or maximum likelihood estimation \cite{beck1977parameter, golightly2008bayesian, myung2003tutorial, pavliotis2007parameter}. Parameters can also be estimated by using the Lagrangian coherent structures \cite{maclean2017coherent} or being treated as the augmented state variables in data assimilation \cite{slivinski2017assimilating}. Here, in light of the mathematically tractable framework, an efficient iteration algorithm is developed that not only provides a systematic estimation of the model parameters but allows the assessment of the uncertainty in a hierarchy of tasks from examining the error in the estimated parameters to its influence on the model forecast skill and the accuracy of data assimilation.
Second, developing cheap approximate filtering strategies is always an alternative to enhancing the data assimilation efficiency \cite{law2012evaluating, cohn1996approximate}. The Lagrangian data assimilation framework allows rigorous mathematical derivations to reach several practical approximate filters, including simplifying the posterior covariance matrix and a randomized selection of observations. Interpreting the uncertainties introduced by these approximations becomes essential for generalizing the strategies to more complicated problems. It is also of practical importance to assess the data assimilation skill of these approximate filters in different scenarios, including recovering both the incompressible and the compressible underlying flows \cite{majda2003introduction, vallis2016geophysical}.
Finally, the data assimilation skill using the analytically solvable filtering formulae with LSMs is compared with that of the ensemble Kalman-Bucy filter \cite{bergemann2012ensemble} with the fully nonlinear forecast model in recovering the highly nonlinear and non-Gaussian geophysical flow fields. Such a study advances the understanding of the error and uncertainty in Lagrangian data assimilation using LSMs for solving many practical problems.

The rest of the paper is organized as follows. The mathematical framework of the nonlinear Lagrangian data assimilation using LSMs is presented {in  Section \ref{Sec:Model} }. The efficient parameter estimation algorithm and its error assessment are included in Section \ref{Sec:ParameterEstimation}. The development of the approximate filtering strategies and the associated UQ are shown in Section \ref{Sec:ReducedFilters}. The comparison of the filtering skill by using linear and nonlinear forecast models in Lagrangian data assimilation is presented in Section \ref{Sec:ComparisonFilters}. The paper is concluded in Section \ref{Sec:Conclusion}.

\section{Mathematical Framework of Nonlinear Lagrangian Data Assimilation with LSMs}\label{Sec:Model}
\subsection{Governing equations of the tracer motion and the underlying flow field}\label{Subsec:GoverningEquations}
Consider a two-dimensional double periodic domain $[0,2\pi)^2$. Assume there are $L$ Lagrangian tracers, whose displacements are denoted by $\mathbf{x}_\ell = (x_\ell, y_\ell)^\mathtt{T}$ with $\ell = 1,\ldots,L$. The governing equation of each $\mathbf{x}_\ell$ satisfies Newton's law
\begin{equation}\label{eq:obs}
\frac{\d \mathbf{x}_\ell(t)}{\d t} = \mathbf{v}(\mathbf{x}_\ell(t), t) + \boldsymbol{\sigma}_\mathbf{x} \dot{\mathbf{W}}_\ell(t),
\end{equation}
where the tracers are treated as massless particles, and therefore the velocity of each tracer $\mathbf{v}(\mathbf{x}_\ell(t), t)$ is the same as that of the {underlying flow field}. The additional white noise forcing with ${\mathbf{W}}_\ell(t)$ being the Brownian motion represents the uncertainty of the tracer movement, accounting for the contribution from the unresolved scales of the flow field. The noise strength $\boldsymbol{\sigma}_\mathbf{x}$ is a $2\times 2$ diagonal matrix with the diagonal entries being $\sigma_x$. Different tracers are assumed to have the same noise strength.
Next, the velocity of the underlying flow field is modeled as
\begin{equation}\label{eq:vel}
\mathbf{v}(\mathbf{x}, {t}) = \mathbf{w}(t) + \sum_{\mathbf{k} \in \mathcal{K}, \alpha \in \mathcal{A}} \widehat{v}_{\mathbf{k},\alpha}(t) {g}_\mathbf{k}(\mathbf{x})\mathbf{r}_{\mathbf{k},\alpha},
\end{equation}
where  $\mathbf{w}(t)$ is the time-varying background velocity field while the remaining terms are the spectral representation of the fluctuation part. The set $\mathcal{K}$ is a finite subset of $\mathbb{Z}^2$ that contains the wavenumbers $\mathbf{k}=(k_x,k_y)$ and the index set $\mathcal{A}$ includes different types of modes associated with the same $\mathbf{k}$. Here, $g_\mathbf{k}(\mathbf{x})$ is the basis function, $\widehat{v}_{\mathbf{k},\alpha}(t)$ is the associated time series, and $\mathbf{r}_{\mathbf{k},\alpha}$ is the eigenvector \cite{chen2015noisy}. Once the family of $g_\mathbf{k}$ is specified, $\mathbf{r}_{\mathbf{k},\alpha}$ is uniquely determined by the underlying physical model. Throughout this paper, the Fourier basis is used, and therefore
\begin{equation}\label{eq:Fourier}
g_\mathbf{k} = \exp(i\mathbf{k}\cdot\mathbf{x}),
\end{equation}
which, together with the expression of the velocity field in \eqref{eq:vel}, implies the governing equation of the tracer movement \eqref{eq:obs} is highly nonlinear.
The following LSMs are utilized to describe the time evolution of the background mean flow $\mathbf{w}(t)$ and the Fourier coefficients  $\widehat{v}_{{\bf k},\alpha}(t)$,
\begin{subequations}\label{eq:uhat}
	\begin{align}
	\frac{\d\mathbf{w}(t)}{\d t} &= -D_0\mathbf{w}(t)+\Omega_0\mathbf{w}(t)+{\bf f}_0+{\boldsymbol\sigma}_0
	\dot{\bf W}_0(t),\label{eq:uhat_w}\\
	\frac{\d\widehat{v}_{{\bf k},\alpha}(t)}{\d t} &= (-d_{{\bf k},\alpha}+i\omega_{{\bf k},\alpha})\widehat{v}_{{\bf k},\alpha}(t)+{f}_{{\bf k},\alpha} +\sigma_{{\bf k},\alpha}
	\dot{W}_{{\bf k},\alpha}(t).\label{eq:uhat_v}
	\end{align}
\end{subequations}
In \eqref{eq:uhat}, $D_0$ and $\Omega_0$ are $2\times 2$ diagonal and anti-diagonal matrices, respectively, ${\bf f}_0$ is a $2\times 1$ vector, ${\boldsymbol\sigma}_0$ is a $2\times 2$ matrix, and $\dot{\bf W}_0(t)$ is a $2\times 1$ white noise. On the other hand, $d_{{\bf k},\alpha}$, $\omega_{{\bf k},\alpha}$, ${f}_{{\bf k},\alpha}$ and $\sigma_{{\bf k},\alpha}$ are all constant scalars and ${W}_{{\bf k},\alpha}(t)$ is a scalar white noise. Among these coefficients and white noises, ${f}_{{\bf k},\alpha}$ and $\dot{W}_{{\bf k},\alpha}(t)$ are complex-valued while all the others are real. The matrix  $D_0$ is positive definite and the scalar $d_{{\bf k},\alpha}$ is a positive number. With the negative signs in front of them, they represent the damping effects. The terms $\Omega_0$ and $\omega_{{\bf k},\alpha}$ introduce wave oscillations and flow rotations, respectively. By choosing suitable coefficients and white noises, each variable $\widehat{v}_{{\bf k},\alpha}(t)$ has its corresponding complex conjugate.

It is worthwhile to remark that if the exact nonlinear system of $\mathbf{v}(\mathbf{x}, {t})$ is known, then each LSM can be calibrated by matching its mean, variance, and decorrelation time with those of the corresponding Fourier mode from the original nonlinear system. However, the exact dynamics of $\mathbf{v}(\mathbf{x}, {t})$ may not be fully available in  practice. Thus, the calibration of these LSMs has to rely on the observed tracer trajectories. The details will be presented in Section \ref{Subsec:ParameterEstimation}.

\subsection{Lagrangian data assimilation}
Lagrangian data assimilation aims to recover the underlying flow field given the observations of tracer trajectories. It can be seen from \eqref{eq:obs}--\eqref{eq:uhat} that once the Fourier coefficients are recovered, the underlying flow field is uniquely determined. Denote by $\mathbf{U}$ a column vector that collects  the  coefficients of the underlying flow fields $\mathbf{w}$ and $\widehat{v}_{{\bf k},\alpha}$ in \eqref{eq:uhat}, which has a dimension {{$K := \vert\mathcal{A}\vert\cdot\vert\mathcal{K}\vert+2$, where $\vert\mathcal{A}\vert$ and $\vert\mathcal{K}\vert$ stand for the number of elements in the sets $\mathcal{A}$ and $\mathbf{K}$, respectively. }} Further denote by $\mathbf{X}$ the collection of the tracer displacement $\mathbf{x}_\ell$ for $\ell=1,\ldots,L$, which has a dimension $2L$. Then the coupled observation-flow system \eqref{eq:obs}--\eqref{eq:uhat} can be written in an abstract form,
\begin{subequations}\label{eq:cgns}
	\begin{align}
	\frac{\d \mathbf{X}(t)}{\d t} &= \mathbf{A}(\mathbf{X}, t) \mathbf{U}(t) + \sigma_x \dot{\mathbf{W}}_\mathbf{X}(t),\label{eq:cgns_X}\\
	\frac{\d \mathbf{U}(t)}{\d t} &=  \mathbf{F}_\mathbf{U} + \boldsymbol{\Lambda} \mathbf{U}(t)  + \boldsymbol{\Sigma}_\mathbf{U} \dot{\mathbf{W}}_\mathbf{U}(t),\label{eq:cgns_U}
	\end{align}
\end{subequations}
where $\mathbf{A}(\mathbf{X}, t)$ contains all the Fourier bases and is, therefore, a highly nonlinear function of $\mathbf{X}$. On the other hand, since the Fourier coefficients $\mathbf{w}(t)$ and $\widehat{v}_{\mathbf{k},\alpha}(t)$ appear linearly by themselves  in \eqref{eq:vel}, there is no nonlinear self interactions of $\mathbf{U}$ in \eqref{eq:cgns}. The conditional linearity of $\mathbf{U}$ in the observational process \eqref{eq:cgns_X} and the LSMs of the flow field \eqref{eq:cgns_U} allow the development of closed analytic formulae for computing the conditional distribution $p(\mathbf{U}(t)|\mathbf{X}(s\leq t))$, where $\mathbf{X}(s\leq t)$ is a given observed realization of the tracer trajectories up to the current time instant $t$. This is also known as the filtering posterior distribution of data assimilation.
\begin{proposition}[Posterior distribution of Lagrangian data assimilation]
	Given one realization of the tracer trajectories $\mathbf{X}(s\leq t)$, the filtering posterior distribution $p(\mathbf{U}(t)|\mathbf{X}(s\leq t))$ of Lagrangian data assimilation \eqref{eq:cgns} is conditionally Gaussian,
	where the time evolutions of the conditional mean $\boldsymbol\mu$ and the conditional covariance $\bf R$ are given by \cite{liptser2013statistics, chen2018conditional}
	\begin{subequations}\label{eq:filter}
		\begin{align}
		\frac{\d\boldsymbol{\mu}}{\d t} &= \left(\mathbf{F}_\mathbf{U} + \boldsymbol{\Lambda} \boldsymbol{\mu}\right)  + \sigma_x^{-2}\mathbf{R}\mathbf{A}^\ast\left(\frac{\d \mathbf{X}}{\d t} - \mathbf{A}\boldsymbol{\mu} \right),\label{eq:filter_mu}\\
		\frac{\d\mathbf{R}}{\d t} &= \boldsymbol{\Lambda}\mathbf{R} + \mathbf{R}\boldsymbol{\Lambda}^\ast + \boldsymbol{\Sigma}_\mathbf{U}\boldsymbol{\Sigma}_\mathbf{U}^\ast - \sigma_x^{-2}\mathbf{R}\mathbf{A}^\ast\mathbf{A}\mathbf{R},\label{eq:filter_R}
		\end{align}
	\end{subequations}
	with $\cdot^*$ being the complex conjugate transpose.
\end{proposition}

\subsection{Assessments of the error and uncertainty}

Both path-wise and information measurements are adopted to assess the error and uncertainty. The former is simple and appropriate to characterize the error in the posterior mean time series, while it is more natural to adopt the latter for evaluating the total uncertainty described by the posterior distribution.

\subsubsection{Path-wise measurements}
The path-wise measurements used here are the normalized root-mean-square error (RMSE) and the pattern correlation (Corr) between the posterior mean states and the truth. The normalized RMSE is the RMSE divided by the standard deviation of the true signal, which gives a non-dimensional quantity. For the conciseness of presentation, the RMSE below always stands for the normalized one.
\begin{definition}[Path-wise measurements]
	The RMSE and Corr are defined as follows \cite{hyndman2006another}:
	\begin{equation}\label{SkillScores}
	\begin{split}
	\mbox{RMSE} &=  \frac{1}{\mbox{std}(u^{ref})}\left(\sqrt{\frac{\sum_{i=1}^I(u^{DA}_{i}-u^{ref}_{i})^2}{I}}\right),\\
	\mbox{Corr} &= \frac{\sum_{i=1}^I(u^{DA}_{i}-\bar{u}^{DA})(u^{ref}_i-\bar{u}^{ref})}{\sqrt{\sum_{i=1}^I(u^{DA}_{i}-\bar{u}^{DA})^2}\sqrt{\sum_{i=1}^I(u^{ref}_{i}-\bar{u}^{ref})^2}},
	\end{split}
	\end{equation}
	where $u^{DA}_{i}$ and $u^{ref}_{i}$ are the posterior mean estimate and the truth of the state variable $u$, respectively, at a single point $i$. The value $I$ is the total number of points for computing these skill scores. Depending on the context, $I$ can be the total number of points in a time series, the total number of the spatial grid points at a fixed time, or the total number of the points in both time and space. The averages of the estimate and the true time series are denoted by $\bar{u}^{DA}$ and $\bar{u}^{ref}$ while std$(u^{ref})$ is the standard deviation of the truth.  A better state estimation corresponds to a smaller RMSE and a larger Corr.
\end{definition}
\subsubsection{Information measurement}
The total uncertainty in the Lagrangian data assimilation can be quantified by computing the difference between the posterior and the prior distributions. Here, the prior distribution is the statistical equilibrium state of the forecast model, which does not involve additional information from observations.
The relative entropy, which is an information criterion, can be utilized to assess the information gain in the posterior distribution $p(t)$ compared with the prior one $p_{eq}$ \cite{majda2005information, kleeman2011information, delsole2005predictability, giannakis2012quantifying, kleeman2002measuring, branicki2013non, majda2002mathematical, chen2014information}.
\begin{definition}[Relative entropy]
	The relative entropy is defined as follows:
	\begin{equation}\label{Relative_Entropy}
	\mathcal{P}(p(t),p_{eq}) = \int p(t)\log\left(\frac{p(t)}{p_{eq}}\right),
	\end{equation}
	which is also known as  Kullback-Leibler divergence  or information divergence \cite{kullback1951information, kullback1987letter}.
	One practical setup in many applications arises when both the distributions involve only the mean and covariance so that
	$p(t)\sim\mathcal{N}(\bar{\mathbf{u}}, \mathbf{R})$ and $p_{eq}\sim\mathcal{N}(\bar{\mathbf{u}}_{eq}, \mathbf{R}_{eq})$
	are Gaussian distributions. In such a case, $\mathcal{P}(p(t),p_{eq})$ has an explicit formula
	\begin{equation}\label{Signal_Dispersion}
	\mathcal{P}(p(t),p_{eq}) = \left[\frac{1}{2}(\bar{\mathbf{u}}-\bar{\mathbf{u}}_{eq})^*\mathbf{R}_{eq}^{-1}(\bar{\mathbf{u}}-\bar{\mathbf{u}}_{eq})\right] + \left[-\frac{1}{2}\log\det(\mathbf{R}\mathbf{R}_{eq}^{-1}) + \frac{1}{2}(\mbox{tr}(\mathbf{R}\mathbf{R}_{eq}^{-1})-K)\right],
	\end{equation}
	where  `det' is the determinant of a matrix. In \eqref{Signal_Dispersion}, the first term in brackets is called `signal', reflecting the information difference in the mean but weighted by the inverse of the equilibrium variance, $\mathbf{R}_{eq}$, whereas the second term in brackets, called `dispersion', involves only the information distance regarding the covariance ratio, $\mathbf{R}\mathbf{R}_{eq}^{-1}$.
\end{definition}
Despite the lack of symmetry, the relative entropy \eqref{Relative_Entropy} has two attractive features. First, $\mathcal{P}(p(t),p_{eq}) \geq 0$ with equality if and only if $p(t)=p_{eq}$. Second, $\mathcal{P}(p(t),p_{eq})$ is invariant under general nonlinear changes of variables.
The signal and dispersion terms in \eqref{Signal_Dispersion} are individually invariant under any (linear) change of variables which maps Gaussian distributions to Gaussians.
In the Lagrangian data assimilation framework \eqref{eq:cgns}, both the prior and posterior distributions are Gaussian. Thus, the formula in \eqref{Signal_Dispersion} will be used to quantify the uncertainty reduction. Since $\mathcal{P}(p(t),p_{eq})$ is time-dependent, the skill score shown below is its time average.

\subsubsection{Relative error in the estimated parameters}
Another useful measurement, especially for quantifying the error in the estimated parameters, is the relative error. Denote by $\boldsymbol\theta^{ref}$ and $\boldsymbol\theta^{est}$ the collections of the true and the estimated parameters.
\begin{definition}[Relative error] The relative error is defined as follows:
	\begin{equation}\label{eq:relative_error}
	\mathcal{E}=\frac{\|\boldsymbol\theta^{est}-\boldsymbol\theta^{ref}\|_{L^2}}{\|\boldsymbol\theta^{ref}\|_{L^2}},
	\end{equation}
	where the standard $L^2$ norm is used in both the denominator and numerator.
\end{definition}
\section{Parameter Estimation and its Error Assessment}\label{Sec:ParameterEstimation}
\subsection{An iterative parameter estimation algorithm}\label{Subsec:ParameterEstimation}
As was mentioned at the end of Section \ref{Subsec:GoverningEquations}, the parameters in the LSMs that characterize the underlying flow field are expected to be estimated utilizing the observed tracer trajectories. The following definitions and  propositions are useful in  parameter estimation.
\begin{definition}[Autocorrelation  and decorrelation time]
	Autocorrelation is the correlation of a signal with a delayed copy of itself as a function of delay \cite{gardiner2009stochastic}. The autocorrelation function (ACF) measures the overall memory of a chaotic system. For a zero mean and stationary time series $u$, the ACF can be calculated as
	\begin{equation}\label{Definition_ACF}
	\mbox{ACF}(t) = \lim_{T\to\infty}\frac{1}{T}\int_0^T\frac{u(t+t')u^*(t')}{\mbox{Var}(u)}dt'.
	\end{equation}
	where $t$ is the delay. The decorrelation time $\tau$ is given by
	\begin{equation}\label{Definition_DecorrTime}
	\tau = \int_{0}^{\infty}\mbox{ACF}(t)dt.
	\end{equation}
\end{definition}
\begin{proposition}[Backward smoother and sampling of the unobserved time series]\label{Prop:Sampling}
	Given one realization of the tracer trajectories $\mathbf{X}(t)$ for $t\in[0,T]$, the smoother estimate $p(\mathbf{U}(t)|\mathbf{X}(s), s\in[0,T])\sim\mathcal{N}(\boldsymbol\mu_\mathbf{s}(t),\mathbf{R}_\mathbf{s}(t))$ of the coupled system is also Gaussian,
	where the conditional mean $\boldsymbol\mu_\mathbf{s}(t)$ and conditional covariance $\mathbf{R}_\mathbf{s}(t)$ of the smoother satisfy the following backward equations\begin{subequations}\label{Smoother_Main}
		\begin{align}
		\frac{\overleftarrow{\d \boldsymbol{\mu}_\mathbf{s}}}{\d t} &=  -\mathbf{F}_\mathbf{U} - \boldsymbol\Lambda\boldsymbol{\mu}_\mathbf{s}  + (\boldsymbol{\Sigma}_\mathbf{U}\boldsymbol{\Sigma}_\mathbf{U}^*)\mathbf{R}^{-1}(\boldsymbol\mu - \boldsymbol{\mu}_\mathbf{s}),\label{Smoother_Main_mu}\\
		\frac{\overleftarrow{\d \mathbf{R}_\mathbf{s}}}{\d t} &= - (\boldsymbol\Lambda + (\boldsymbol{\Sigma}_\mathbf{U}\boldsymbol{\Sigma}_\mathbf{U}^*) \mathbf{R}^{-1})\mathbf{R}_\mathbf{s} - \mathbf{R}_\mathbf{s}(\boldsymbol\Lambda^* + (\boldsymbol{\Sigma}_\mathbf{U}\boldsymbol{\Sigma}_\mathbf{U}^*)\mathbf{R})  + \boldsymbol{\Sigma}_\mathbf{U}\boldsymbol{\Sigma}_\mathbf{U}^* ,\label{Smoother_Main_R}
		\end{align}
	\end{subequations}
	with $\boldsymbol\mu$ and $\mathbf{R}$ being given by \eqref{eq:filter}. The notation $\overleftarrow{\d \cdot}/\d t$ corresponds to the negative of the usual derivative, which means that the system \eqref{Smoother_Main} is solved backward over $[0,T]$ with the starting value of the nonlinear smoother being the same as the filter estimate $(\boldsymbol\mu_\mathbf{s}(T), \mathbf{R}_\mathbf{s}(T)) = (\boldsymbol\mu(T), \mathbf{R}(T))$.
	Based on the smoother estimate, an optimal backward sampling of the trajectories associated with the unobserved variable $\mathbf{U}$ satisfies the following explicit formula \cite{chen2020efficient},
	\begin{equation}\label{Sampling_Main}
	\frac{\overleftarrow{\d \mathbf{U}}}{\d t} = \frac{\overleftarrow{\d \boldsymbol\mu_\mathbf{s}}}{\d t} - \big(\boldsymbol\Lambda + (\boldsymbol{\Sigma}_\mathbf{U}\boldsymbol{\Sigma}_\mathbf{U}^*)\mathbf{R}^{-1}\big)(\mathbf{U} - \boldsymbol\mu_\mathbf{s}) + \boldsymbol{\Sigma}_\mathbf{U}\dot{\mathbf{W}}_{\mathbf{U}}(t).
	\end{equation}
\end{proposition}
\begin{proposition}[Parameter estimation by matching statistics]
	The four parameters $f_{\bf k,\alpha}$, $d_{\bf k,\alpha}$, $\omega_{\bf k,\alpha}$ and $\sigma_{\bf k,\alpha}$ in the complex-valued linear stochastic model \eqref{eq:uhat_v} can be uniquely determined by utilizing the four statistics of $\hat{v}_{\bf k,\alpha}$: the mean $m_{\bf k,\alpha}$, the variance $E_{\bf k,\alpha}$, and the real and imaginary parts of the decorrelation time $T_{\bf k,\alpha}- i\theta_{\bf k,\alpha}$ \cite{majda2012filtering},
	\begin{equation}\label{eq:matching_statistics}
	f_{\bf k,\alpha}=\frac{(T_{\bf k,\alpha}-i\theta_{\bf k,\alpha})m_{\bf k,\alpha}}{T_{\bf k,\alpha}^2+\theta_{\bf k,\alpha}^2},~
	d_{\bf k,\alpha}=\frac{T_{\bf k,\alpha}}{T_{\bf k,\alpha}^2+\theta_{\bf k,\alpha}^2},~
	\omega_{\bf k,\alpha}=\frac{\theta_{\bf k,\alpha}}{T_{\bf k,\alpha}^2+\theta^2_{\bf k,\alpha}},~
	\sigma_{\bf k,\alpha}=\sqrt{\frac{2E_{\bf k,\alpha}T_{\bf k,\alpha}}{T_{\bf k,\alpha}^2+\theta_{\bf k,\alpha}^2}}.
	\end{equation}
	On the other hand, the parameters in the two-dimensional real-valued large-scale mean flow equations \eqref{eq:uhat_w} can be easily estimated by utilizing a simple linear regression.
\end{proposition}

With these propositions in hand, the iterative algorithm of estimating the parameters in the linear stochastic models \eqref{eq:cgns} is presented in Algorithm \ref{Alg:ParameterEstimation}.
\begin{remark}[Computing the decorrelation time in practice]
	In practice, only a finite length of the time series is available to numerically compute the ACF \eqref{Definition_ACF} and the decorrelation time \eqref{Definition_DecorrTime}. As a result, small fluctuations always exist at the tail of the ACF, and the accumulated effect of these fluctuations can often introduce a significant error in computing the decorrelation time. Implementing an accurate finite-time truncation of the ACF curve is difficult in practice without seeing the profile of the ACF, which is the case in the iterative parameter estimation algorithm. Nevertheless, the ACF of the LSM has an analytic ansatz, where the cross-correlation function between the real and
	imaginary parts equals $\exp(-c_1t)\sin(c_2t)$, with $c_1>0$ and $c_2$ being the damping and oscillation rates of the time series. Thus,  the numerically computed ACF is used to determine the two constants in the above ansatz directly related to the $T_{\bf k,\alpha}$ and $\theta_{\bf k,\alpha}$ in the decorrelation time.
\end{remark}
\begin{remark}[The necessity of the backward sampling]
	Although the backward sampling in Proposition \ref{Prop:Sampling} slightly increases the computational cost, it is essential to unbiasedly compute the statistics of each Fourier mode. A cruder approximation in obtaining the trajectories of the unobserved variables during the parameter estimation procedure is to use the filter mean time series. However, the filtering solution is generally less accurate than smoothing, as the former takes into account the observational information only in the past. In addition, the filter or smoother mean time series underestimates the variability of the recovered variable. Including the information in the posterior covariance and the temporal dependence, as in the sampled trajectories, is crucial to obtain unbiased statistics.
\end{remark}

\begin{algorithm}
	\caption{An iterative parameter estimation algorithm for linear stochastic models}
	\label{Alg:ParameterEstimation}
	\begin{algorithmic}
		\State{Given $\mathcal{K},\mathcal{A},\sigma_x$ and the trajectories of $L$ tracers. Set up an error threshold $\epsilon$.}
		\State{Start with an initial guess of the parameters, denoted by $\boldsymbol\theta_0$. Set $\boldsymbol\theta_1=\boldsymbol\theta_0-2\epsilon$ and $n=1$.}
		\While{$\|\boldsymbol\theta_{n}-\boldsymbol\theta_{n-1}|\geq\epsilon$}
        \State{{Compute the filter and smoother solutions \eqref{eq:filter} and \eqref{Smoother_Main}.}}
		\State{{{Based on the filter covariance $\mathbf{R}$ and the smoother mean $\boldsymbol\mu_\mathbf{s}$}}, apply \eqref{Sampling_Main} to obtain a sampled trajectory of $\mathbf{U}$.}
		\State{For each $\hat{v}_{\bf k,\alpha}$, use the corresponding component of $\mathbf{U}$ to compute the mean, variance and decorrelation time {using the formulae in \eqref{eq:matching_statistics}}. }
		\State{Use \eqref{eq:matching_statistics} to update the parameters in the equation of $\hat{v}_{\bf k,\alpha}$ based on these statistics. }
		\State{Use a linear regression to determine the parameters in $\mathbf{w}(t)$. }
		\State{Denote the new parameters by $\boldsymbol\theta_{n+1}$ and let $n := n+1$.}
		%\STATE{$n:=n+1$.}
		\EndWhile\\
		\Return $\boldsymbol\theta_{n+1}$
	\end{algorithmic}
\end{algorithm}

\subsection{Assessing the uncertainty in the estimated parameters and the associated model}\label{Subsec:UQ_Param}
The following model is utilized to study the skill of the parameter estimation and the associated Lagrangian data assimilation skill. The underlying flow field is given by a set of LSMs \eqref{eq:uhat}, and the same model ansatz is adopted for parameter estimation and data assimilation to eliminate the structural model error. The underlying flow field does not contain the background mean flow, while the fluctuation part includes Fourier wavenumbers ranging from $-5\leq k_x,k_y\leq 5$. The flow field is incompressible, representing the geophysically balanced (GB) features, which means the set of $\mathcal{A}$ has only one element, and therefore, the index $\alpha$ is omitted here. {{In such a case, there are in total $K=121$ modes.}} The eigenvector of the incompressible flow is $\mathbf{r}_{\bf k} = (-ik_y, ik_x)/\vert\mathbf{k}\vert^2$.  The true parameters are as follows: $f_{\bf k}=0$,  $\omega_{\bf k}=0$,  $d_{\bf k}=d+\nu |{\bf k}|^2$, with $d=0.3,\nu=0.05$. Let $E_0=1$, $\alpha=3$, and $k_0=2$. An energy spectrum is assigned for the flow field,
\begin{equation}
E_\mathbf{k} =
\begin{cases}
\vert\mathbf{k}\vert E_0, & \vert\mathbf{k}\vert\leq k_0, \\
k_0 E_0 \left|\frac{\vert\mathbf{k}\vert}{k_0}\right|^{-\alpha}, & \vert\mathbf{k}\vert > k_0,
\end{cases}
\end{equation}
which is used to compute the stochastic forcing coefficients via the relationship \cite{majda2012filtering}
\begin{equation}\label{eq:spectrum}
E_{\mathbf{k}} = \frac{1}{2}\left(\frac{f_{\mathbf{k}}^2}{d_{\mathbf{k}}^2} + \frac{\sigma_{\mathbf{k}}^2}{2 d_{\mathbf{k}}}\right),
\end{equation}
The noise coefficient in the observational processes is $\sigma_x = 0.25$. A simulation with $400$ time units is utilized, which is about $120$ times the decorrelation time of the slowest decaying mode and is sufficiently long to exclude most of the random error in computing the skill scores. The numerical integration time step is $\Delta{t} = 0.002$ for both the model simulation and data assimilation. The Euler-Maruyama scheme is adopted for numerical integration.

\begin{figure}[htbp]
	\centering
	\hspace*{-0.0cm}\includegraphics[width=1.00\textwidth]{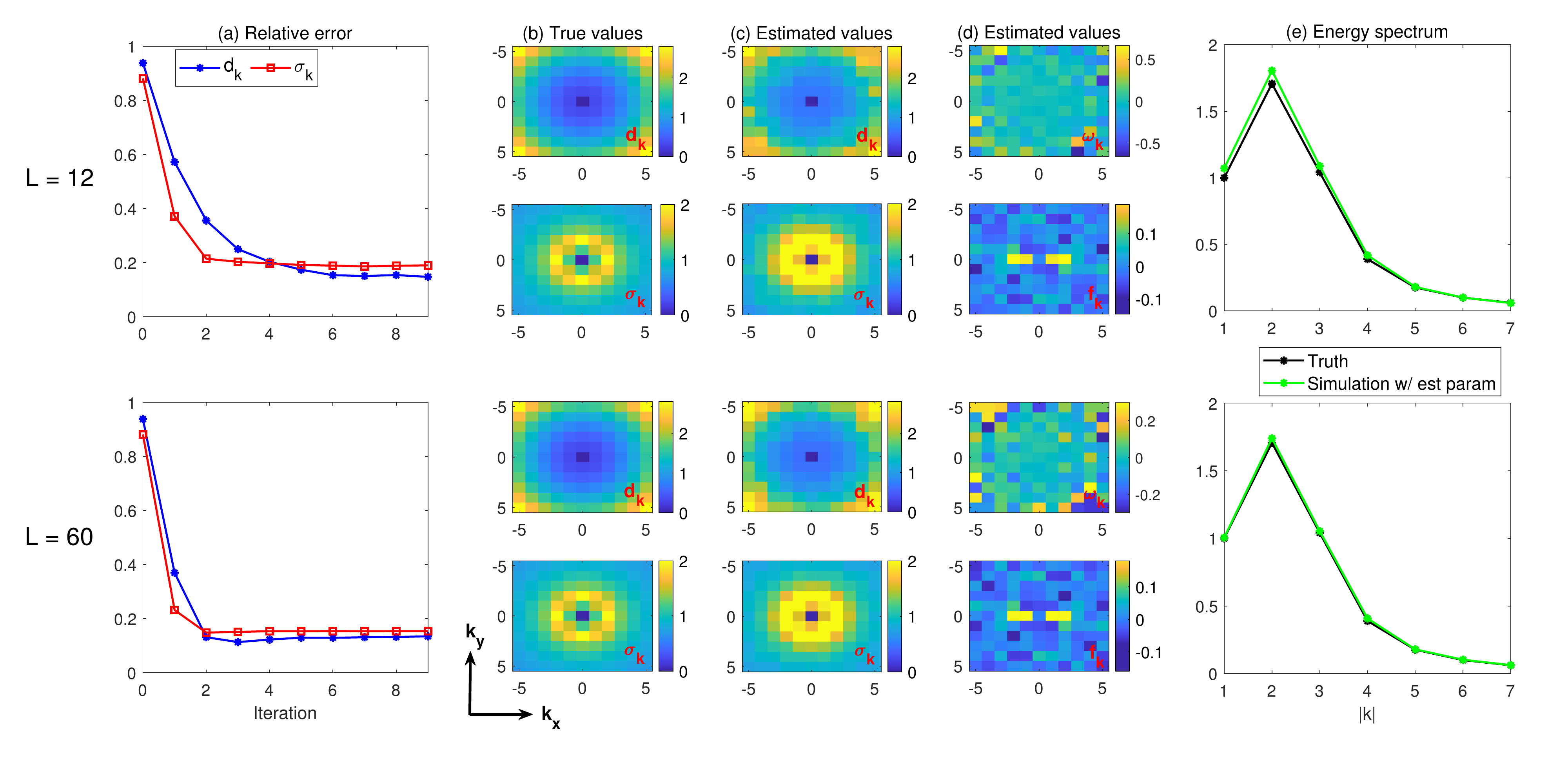}
	\caption{Assessment of the parameter estimation skill. Panel (a): the relative error in the estimated damping and noise coefficients, $d_\mathbf{k}$ and $\sigma_\mathbf{k}$, as a function of iterations. Here, the error is computed by exploiting \eqref{eq:relative_error}, where the $d_\mathbf{k}$ or the $\sigma_\mathbf{k}$ for all the $121$ modes are included into the vectors $\boldsymbol\theta^{ref}$ and $\boldsymbol\theta^{est}$. Panel (b): the true parameter values of $d_\mathbf{k}$ and $\sigma_\mathbf{k}$ for different wavenumbers. Panel (c): the estimated $d_\mathbf{k}$ and $\sigma_\mathbf{k}$. Panel (d): the estimated parameters $\omega_{\mathbf{k}}$ and $f_{\mathbf{k}}$, { where the true values of $\omega_{\bf k}$ and $f_{\bf k}$ are all zeros}. Panel (e): energy spectrum as a function of $\mathbf{k}$. The black and green curves show the spectrums using the formula \eqref{eq:spectrum} with the true parameters and the estimated ones, respectively. }
	\label{fig:gbfig1}
\end{figure}

The assessment of the parameter estimation skill is included in Fig.\ref{fig:gbfig1}. Panel (a) shows that the estimated parameters from  Algorithm \ref{Alg:ParameterEstimation} have a quick convergence within only a few iterations. The error in the estimated parameters is already small with the help of only a few tracers. However, the error does not completely vanish even with a relatively large number of observations $L=60$ due to the noise and the strong nonlinear coupling of different modes in the observational processes. In addition, as is shown in Panels (b)--(d), the estimated parameters of the large-scale modes (i.e., small $|\mathbf{k}|$), which account for most of the energy of the system, are slightly more accurate than those of the medium- and small-scale modes. Despite such errors in the estimated parameters, the energy spectrum computed from \eqref{eq:spectrum} of the system using the true parameters and the estimated ones are very similar, { as small errors in the point-wise estimates can be averaged out in computing the statistics}. In fact, although both $\sigma_\mathbf{k}$ and $d_\mathbf{k}$ are slightly overestimated for the medium- and small-scale modes, their ratio $\sigma_\mathbf{k}^2/(2d_\mathbf{k})$ that contributes to the total energy is accurately recovered. {This is significant as recovering the energy spectrum is often an important task in practice.} The results here justify the skill of the parameter estimation algorithm.

\begin{figure}[htbp]
	\centering
	\hspace*{-0.0cm}\includegraphics[width=1.00\textwidth]{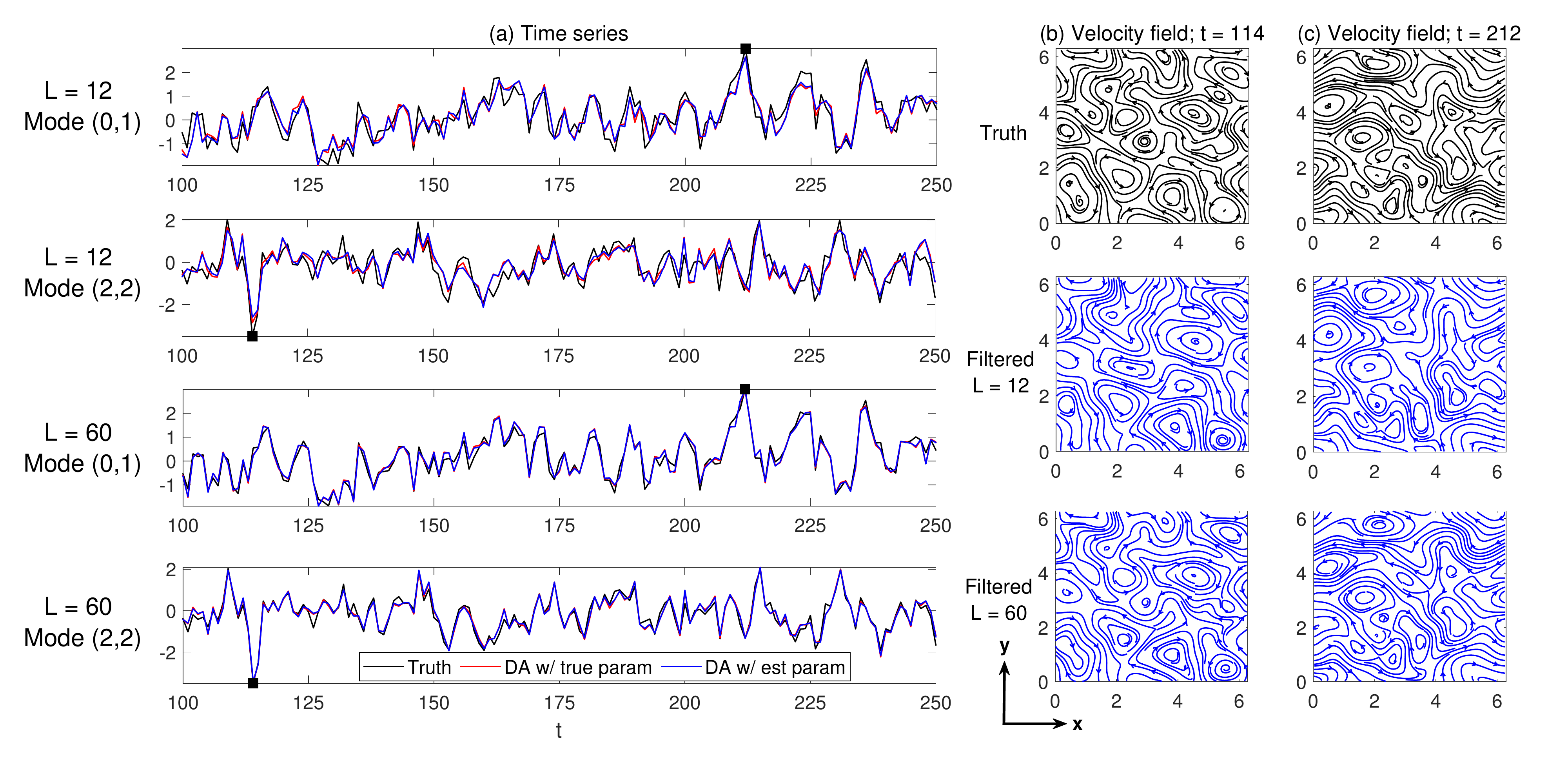}
	\caption{Comparison of the posterior mean estimates. Panel (a): the true signal (black) and the posterior mean time series from the Lagrangian data assimilation using the true parameters (red) and the estimated ones (blue). Panel (b): the true and the filtered velocity fields at time $t=114$, where the filtered ones are computed using the estimated parameters. Panel (c): similar to Panel (b) but at $t=212$. These two time instants are marked by black squares in the time series shown in Panel (a).}
	\label{fig:gbfig2}
\end{figure}

Fig.\ref{fig:gbfig2} compares the posterior mean time estimates. Panel (a) shows the time series of a large-scale mode $(0,1)$ and a medium-scale one $(2,2)$. The two posterior mean time series using the true and estimated parameters are quite close to each other, even with only $L=12$ tracers. The difference between the posterior mean time series and the truth decreases as $L$ goes from $L=12$ to $L=60$, as is expected. Panels (b)--(c) compare the true and the filtered velocity fields, where the filtered ones are computed using the estimated parameters. At these two time instants, a strong signal occurs in mode $(2,2)$ and mode $(0,1)$, respectively. It can be seen that the large-scale features are already captured with even $L=12$ tracers. When $L=60$, improvements are found in recovering many small-scale features.

\begin{figure}[htbp]
	\centering
	\hspace*{-0.0cm}\includegraphics[width=1.00\textwidth]{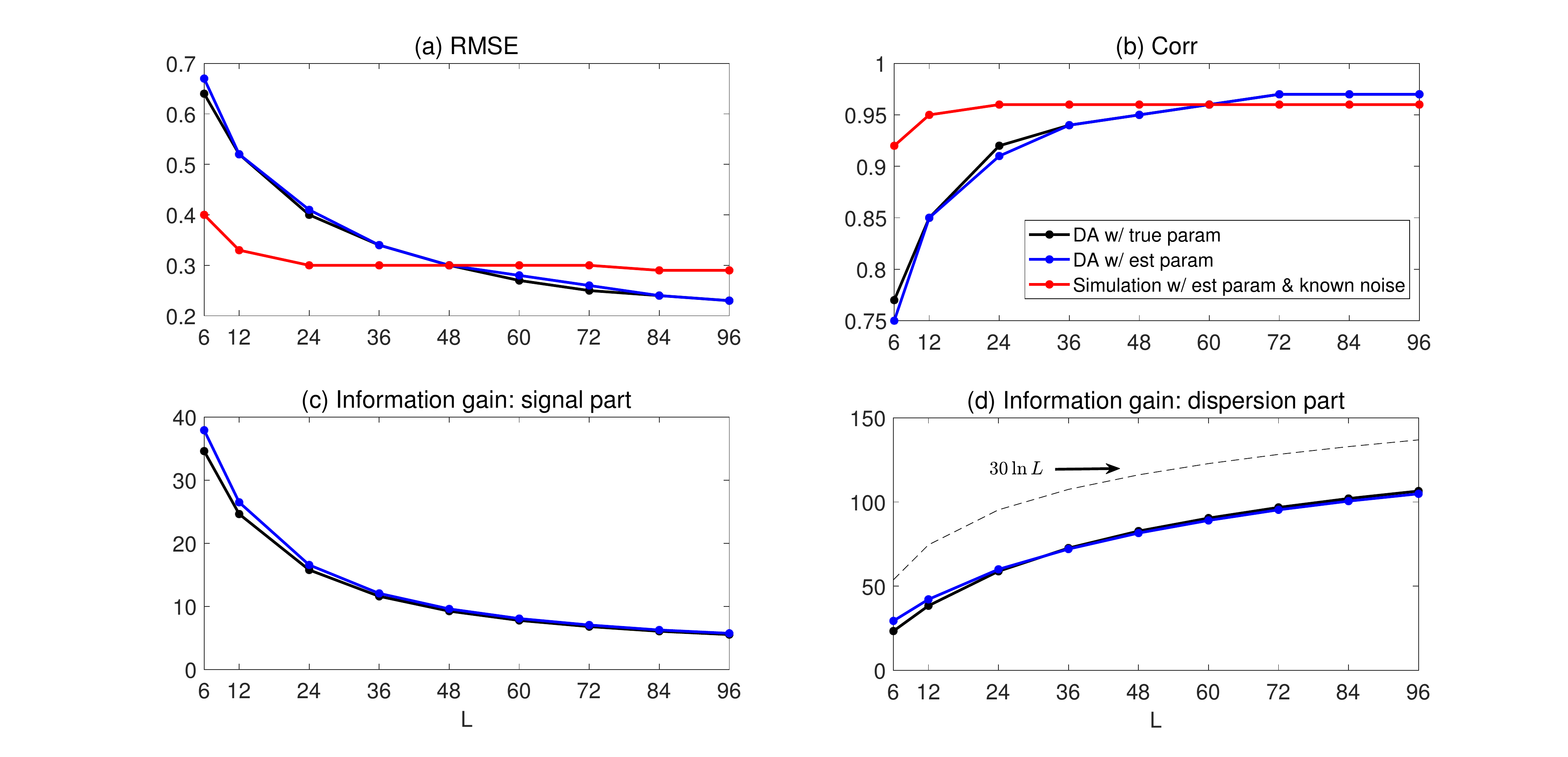}
	\caption{Path-wise and information skill scores as a function of the number of the observed tracers $L$. Panels (a)--(b): RMSE and Corr in the posterior mean time series compared with the truth. The black and blue curves show the results where the posterior time series is computed using the true and estimated parameters, respectively. The red curve shows the skill score from the direct model simulation using the estimated parameters, where the white noise of each Fourier mode is assumed to be the same as that in generating the true signal. Panels (c)--(d): information gain in the posterior distribution related to the prior. The dashed line in Panel (d) represents 30$\ln L$ for comparison. { Note that the number 30 here does not have a special meaning. It is chosen to make the dashed curve different from two solid ones such that the $\ln L$ behavior can be easily seen.}}
	\label{fig:gbfig3}
\end{figure}

Finally, Fig.\ref{fig:gbfig3} illustrates the path-wise and the information skill scores as a function of  $L$. It again confirms that the posterior states using the true and the estimated parameters are almost identical. The path-wise skill scores in Panels (a)--(b) show that the errors monotonically decrease as $L$. Note that the red curve shows the skill score from the direct model simulation using the estimated parameters, where the white noise of each Fourier mode is assumed to be the same as the model that generates the true trajectory. Therefore, the error solely comes from the parameter estimation. The results here show that when $L$ is small, the relatively small error in the direct model simulation using the estimated parameters and the known noise trajectories implies the skill of parameter estimation. On the other hand, when $L$ becomes large, the observational information advances the improvement of the state estimation in the posterior time series, which outweighs the model simulation even with the known random noise sources. Panels (c)--(d) show the information measurements. The signal part converges to a fixed value as the posterior mean time series converges to the truth while the dispersion part increases as a function of $\text{log}(L)$ \cite{chen2014information}. The latter indicates the gain of information in the posterior covariance, which was often ignored when only the path-wise measurements were utilized in quantifying the data assimilation skill.

\section{Approximate Filters and the Associated UQ}\label{Sec:ReducedFilters}
In practice, a large number of spectral modes is needed to characterize the multiscale feature of the underlying complex flow field. When the dimension of the state variables of the flow field $K$ increases, computing the posterior covariance $\mathbf{R}$ in \eqref{eq:filter}, which is of size $K^2$, will become { prohibitively} expensive. Therefore, developing suitable approximate filters that accelerate the computational efficiency of Lagrangian data assimilation is extremely useful. The analytically tractable features of the Lagrangian data assimilation framework with LSMs allow the development of various approximate filters and the quantification of the associated uncertainty. For the convenience of discussion, the background flow $\mathbf{w}(t)$  in \eqref{eq:vel} is temporally omitted in the discussions of this section.

\subsection{Strategies of developing approximate filters}
Although different spectral modes are fully coupled in the observational processes \eqref{eq:cgns_X}, the LSMs of these modes are by design independent in the forecast model \eqref{eq:cgns_U}. Therefore, one natural strategy for developing approximate filters is to assume the posterior covariance in \eqref{eq:filter_R} is a diagonal matrix, where only the $K$ diagonal entries need to be updated. Define $\mathbf{P}=\mathbf{A}^*\mathbf{A}$ in \eqref{eq:filter_R}. The entry of $\mathbf{P}$ that corresponds to the interaction between modes $\mathbf{m}$ and $\mathbf{n}$ can be explicitly written as
\begin{equation}\label{eq:matrixP}
\mathbf{P}_{\mathbf{m},\mathbf{n}} = \sum_{l=1}^L\exp(i(\mathbf{m}-\mathbf{n})\mathbf{x}_l) \mathbf{r}^*_{\mathbf{n}}\mathbf{r}_{\mathbf{m}}.
\end{equation}
{{Note that for the conciseness of notation, the index $\mathbf{m}$ and $\mathbf{n}$ contains not only the modes $\bf k$ but also the index $\alpha$.}} Denote by $r_{\bf k,\alpha}$ and $\mathbf{P}_{\bf k,\alpha}$ the diagonal entry of the covariance matrix and the entry of $\mathbf{P}_{\mathbf{m},\mathbf{n}}$ associated with $\hat{v}_{\bf k,\alpha}$, respectively. Then the evolution equation of $r_{\bf k,\alpha}$ satisfies
\begin{equation}\label{eq:diag_cov}
\frac{\d r_{\bf k,\alpha}}{\d t} = 2d_{\bf k,\alpha} r_{\bf k,\alpha}  + \sigma_{\bf k,\alpha}^2 - \sigma_x^{-2}\mathbf{P}_{\bf k,\alpha}r_{\bf k,\alpha}^2.
\end{equation}
\begin{definition}[Approximate filter I: Diagonal posterior covariance matrix]\label{Prop:ApproximateFilter_I}
	The posterior distribution of such an approximate filter is given by \eqref{eq:filter_mu} and \eqref{eq:diag_cov}.
\end{definition}

It has been shown in \cite{chen2014information} that if the underlying flow is incompressible, then the statistical equilibrium distribution of the tracers' locations converges geometrically fast towards the uniform distribution in $[0,2\pi)^2$. Therefore, the entries of $\mathbf{P}_{\mathbf{m},\mathbf{n}}$ with $\mathbf{m}\neq\mathbf{n}$ become zero in the asymptotic limit according to the mean field theory while those with $\mathbf{m}=\mathbf{n}$ equal to $L$,
\begin{equation}\label{eq:diag_cov2}
\frac{\d r_{\bf k}}{\d t} = 2d_{\bf k} r_{\bf k}  + \sigma_{\bf k}^2 - L\sigma_x^{-2}r_{\bf k}^2,
\end{equation}
where the index $\alpha$ has been omitted since the set $\alpha\in \mathcal{A}$ in \eqref{eq:vel} has only one element. The equilibrium solution of $r_{\bf k,\alpha}$ is given by
\begin{equation}\label{eq:diag_cov3}
r_{\bf k} = \frac{\sigma_{\bf k}^2}{d_{\bf k} + \sqrt{d_{\bf k} + L\sigma_x^{-2}\sigma_{\bf k}^2}}
\end{equation}
Based on Proposition \ref{Prop:ApproximateFilter_I}, an even simpler and cruder approximate filter includes a constant posterior covariance matrix is given as follows.
\begin{definition}[Approximate filter II: Diagonal and constant posterior covariance matrix]
	Only the posterior mean \eqref{eq:filter_mu} needs to be updated in such an approximate filter while the posterior covariance is fixed as in \eqref{eq:diag_cov3}.
\end{definition}

In addition to the update of the covariance matrix, another major source of the computational burden comes from the calculation of the matrix $\mathbf{P}=\mathbf{A}^*\mathbf{A}$, which appears in the governing equations of both the mean and the covariance. Since the size of $\mathbf{A}$ is $L\times K$, the  cost of computing $\mathbf{P}$ grows as a quadratic function of the number of tracers $L$. %In other words, although the accuracy of the Lagrangian data assimilation improves as $L$ becomes larger, the computational cost increases as well.
Nevertheless, the continuous observation of the tracer trajectories allows the development of a randomized selection of observations that can significantly reduce the computational cost of $\mathbf{P}$ and thus the posterior distribution.
\begin{definition}[Approximate filter III: randomized selection of observations]
	At each time instant $t_j = j\Delta{t}$, a random number generator is applied to select $S$ out of in total $L$ tracer trajectories, where $S\ll L$ and $\Delta{t}\ll 1$ is the numerical integration time step. Only the displacements of these $S$ trajectories at $t_j$ and $t_{j-1}$ are used in the Lagrangian data assimilation.
\end{definition}
Clearly, throughout the entire data assimilation process, the cost of computing $\mathbf{P}$ decreases from $O(L^2)$ to $O(S^2)$. The following proposition states that, with a minor  adjustment, the posterior mean estimate using such an approximate filter converges to that using the filter with the full observations when $\Delta{t}\to 0$ for incompressible flows.
\begin{proposition}[Convergence of the  approximate filter with randomly selected observations]\label{theorem:random}Denote by $(\boldsymbol\mu_{\bf I},\mathbf{R}_{\bf I})$ and $(\boldsymbol\mu_{\bf II},\mathbf{R}_{\bf II})$ the posterior estimate from the full filter and the approximate one with randomly selected observations. There are in total $L\gg 1$ tracers while only $L'\leq L$ tracers are used at each time instant in the approximate filter. Assume that all the initial values of these two filters are the same. Further assume that, within a short time interval $[t, t+L\Delta{t}]$, the changes of the displacements of all the tracers are small, e.g., $\vert\mathbf{x}_l(t+L\Delta{t})-\mathbf{x}_l(t)\vert\leq \epsilon$ as $\Delta{t}\to 0$. Then
	a) the posterior covariance of the approximate filter converges  to a constant diagonal matrix with $L'\to\infty$, where the diagonal entry is
	\begin{equation}
	r_{\bf II, k} = \sigma_{\bf k}^2/\sqrt{d_{\bf k} + L'\sigma_x^{-2}\sigma_{\bf k}^2},
	\end{equation}
	and b) after multiplying a rescaling prefactor $\sqrt{L/L'}$ in front of the gain function of the posterior mean equation of the approximate filter,
	\begin{equation}\label{eq:prefactor}
	\frac{\d\boldsymbol{\mu}_{\bf II}}{\d t} = \left(\mathbf{F}_\mathbf{U} + \boldsymbol{\Lambda} \boldsymbol{\mu}_{\bf II}\right)  + \sqrt{\frac{L}{L'}}\sigma_x^{-2}\mathbf{R}_{\bf II}\mathbf{A}^\ast\left(\frac{\d \tilde{\mathbf{X}}}{\d t} - \mathbf{A}\boldsymbol{\mu}_{\bf II} \right),
	\end{equation}
	the difference in the posterior mean estimates between using the approximate filter and using the full filter is bounded at any future time with
	\begin{equation}\label{eq:random_obs_bounds}
	\vert\boldsymbol\mu_{\bf I}(t)-\boldsymbol\mu_{\bf II}(t)\vert\leq C\epsilon.
	\end{equation}
	Here $\tilde{\mathbf{X}}$ contains $L'$ randomly selected tracers.
\end{proposition}

The justification of this proposition is shown in Appendix. Essentially, within a short time interval $[0,L\Delta t]$, the effect of the random selection strategy is equivalent to using a large number of tracers $L'L$. It is also worth highlighting that when $L$ is ample in both the full filter and the Approximate filter III, the error in the observations is essentially averaged out. Therefore, the role of such a prefactor in the Approximate filter III is to assign more weight to the information provided by the observations in the filter solution. Note that despite the requirement of $L'\to\infty$ in the justification of the proposition, a small value of $L'$ already works sufficiently well in practice, as shown in the following numerical simulations.

\subsection{Test models}
Two test models are utilized here to study the skill of the approximate filters in Lagrangian data assimilation. The first test model is the same as that utilized in Section \ref{Subsec:UQ_Param}, which describes an incompressible flow field and contains only the GB modes, except that now $-3\leq k_x,k_y\leq3$ such that the total number of Fourier modes is $K = 49$. The second test model is a linear shallow water system, which characterizes a compressible flow field. The model is included in Appendix. For each wavenumber $\mathbf{k}$ in the shallow water equation, the set $\mathcal{A}$ contains three elements: one GB mode $\alpha:= B$ and a pair of gravity modes $\alpha:=\pm$. Thus, there are in total $K=147$ modes. The GB mode describes the incompressible part of the overall flow field, while the two gravity modes contribute to the compressible part. The damping $d_{{\bf k}, \pm}$ and the forcing $\mathbf{f}_{\bf k, \pm}$ of the gravity modes are the same as those in the GB modes $d_{{\bf k}, B}$ and $\mathbf{f}_{{\bf k}, B}$, { as  described in Section \ref{Subsec:UQ_Param}}. The oscillation frequency is given by $\omega_{{\bf k}, \pm}=\pm\mbox{Ro}^{-1}\sqrt{|k|^2+1}$ where $\mbox{Ro}$ is the Rossby number that describes the ratio of inertial force to Coriolis force. It is set to be $\mbox{Ro} = 1$ here. Such a moderate oscillation represents a tough test case, where the gravity and GB modes of the same wavenumber occur in a comparable time scale. The energy spectrum of the gravity modes has the same profile as that of the GB modes. But the energy associated with each gravity mode is $1/4$ of that associated with the GB mode such that the total energy in the GB part of the flow is twice as much as that in the gravity part, which is a reasonable ratio in many geophysical scenarios. See \cite{majda2003introduction, vallis2017atmospheric, chen2015noisy} for more details.

\subsection{Lagrangian data assimilation using approximate filters}
\begin{figure}[htbp]
	\centering
	\hspace*{-0.0cm}\includegraphics[width=1.00\textwidth]{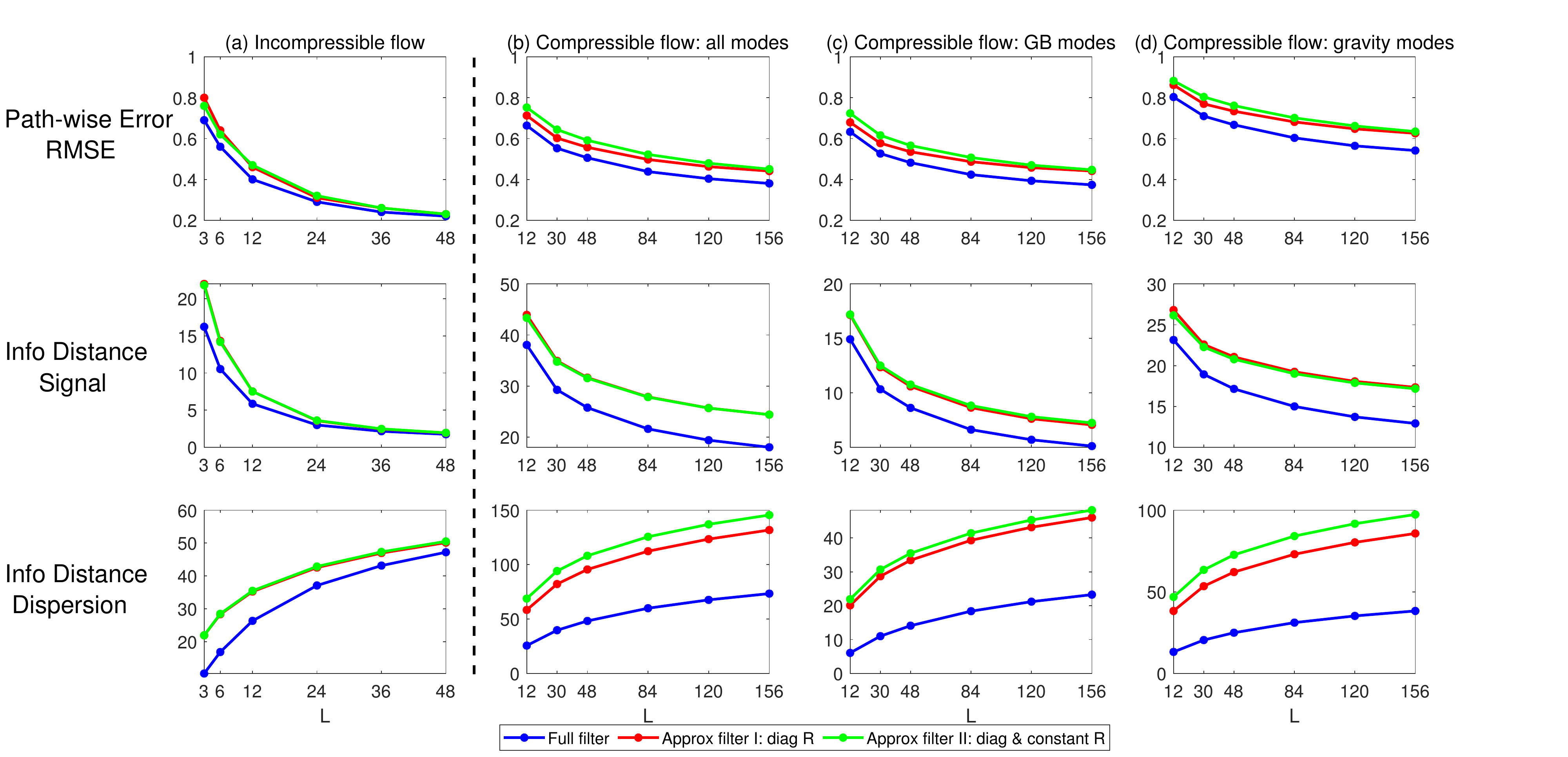}
	\caption{Skill scores using the full filter (blue), the Approximate filter I with a diagonal covariance matrix (red), and the Approximate filter II with a diagonal and constant covariance matrix (green). Column (a) shows the case of the incompressible flow model, while Columns (b)--(d) show that of the compressible shallow water equation. The skill scores in Column (b) are computed based on all the shallow water equation modes, including both the GB and the gravity ones. The skill scores in Column (c) are calculated only based on the filtered GB modes in the shallow water equation. Still, the Lagrangian data assimilation is applied to the entire shallow water equation. Column (d) is similar to Column (c) but for the filtered gravity modes.}
	\label{fig:covfig1}
\end{figure}

Let us start by comparing the filtering skill between the full filter and the Approximate filters I and II, which allows us to understand the additional uncertainty due to the approximation of the posterior covariance matrix. The results are presented in Fig.\ref{fig:covfig1}. Column (a) shows the case of the incompressible flow model, which contains only the GB modes. It is seen that the two approximate filters have almost the same data assimilation skill as the truth in terms of the path-wise RMSE in the posterior mean estimate. The similarity in the posterior mean time series naturally gives a comparable information gain in the signal part. Next, since the approximations are imposed in the posterior covariance, the approximate filters display a more significant dispersion than the full filter. In other words, the estimate is overconfident towards the posterior mean time series. It is worthwhile to highlight that the additional values in the dispersion are the model error rather than the actual uncertainty reduction from Lagrangian data assimilation.

\begin{figure}[htbp]
	\centering
	\hspace*{-0.0cm}\includegraphics[width=1.00\textwidth]{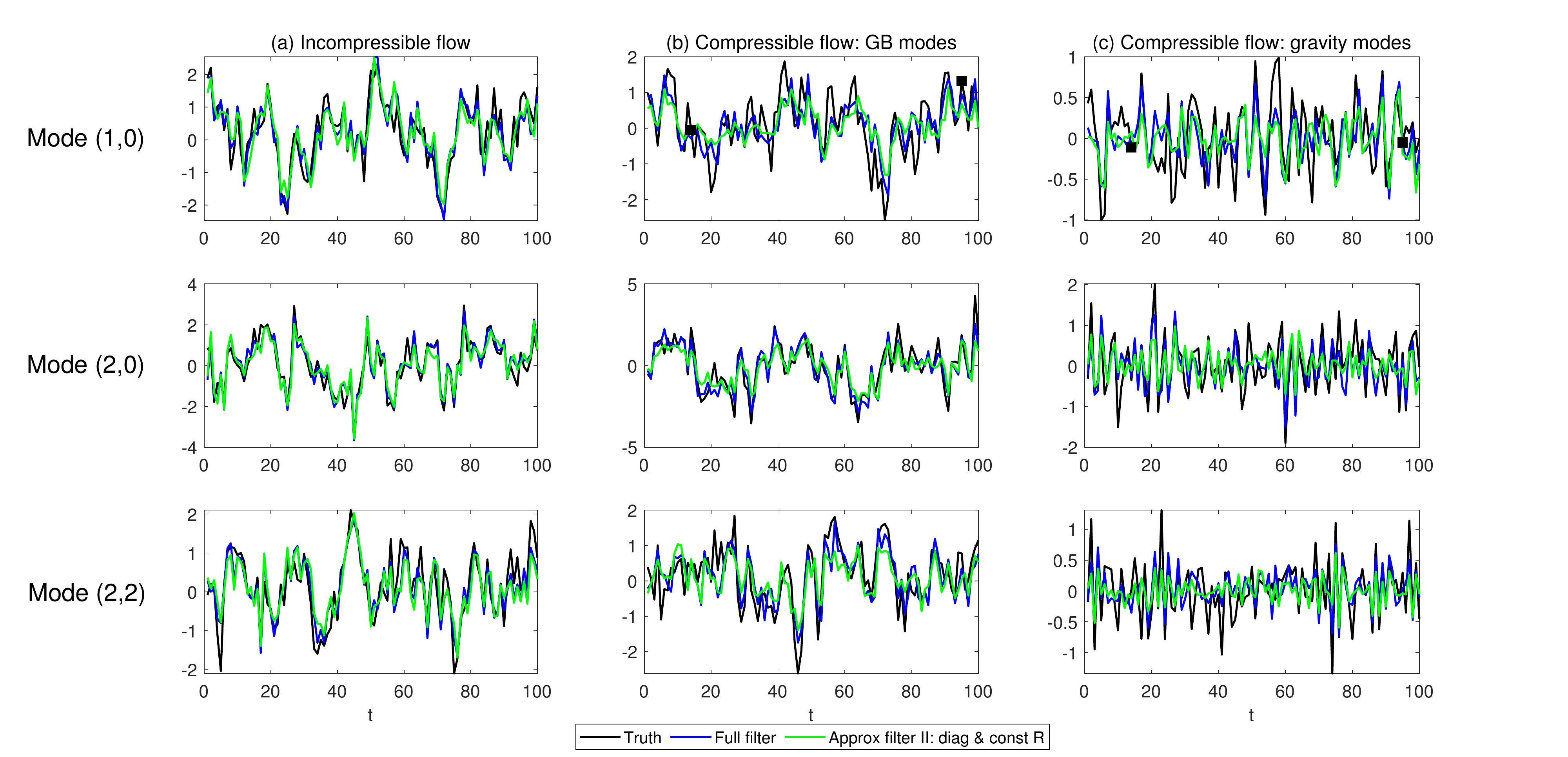}
	\caption{Comparison of the true signal (black) with the posterior mean time series using the full filter (blue) and the Approximate filter II with a diagonal and constant covariance matrix (green) for three modes $\mathbf{k}=(1,0)$, $(2,0)$ and $(2,2)$. Column (a) shows the case of the incompressible flow model, which contains only the GB modes, while Columns (b)--(c) show the results for the GB and the gravity modes of the compressible shallow water equation, respectively.}
	\label{fig:covfig2}
\end{figure}

Columns (b)--(d) of Fig.\ref{fig:covfig1} show the results of the compressible shallow water equation. The difference in the skill scores between the approximate filters and the full filter becomes more pronounced. In particular, the approximate filters are inaccurate in characterizing the uncertainty in the posterior distribution. This is not surprising since the approximate filters based on the asymptotic formulae \eqref{eq:diag_cov2}--\eqref{eq:diag_cov3} are only { mathematically valid for the incompressible case }, where the associated derivations exploit the property of the uniform distribution of the tracers due to the flow incompressibility. Next, Columns (c) and (d) show that it is more challenging to recover the gravity modes than the GB ones, as the structures of the former are more complicated. Because the GB and the gravity modes are mixed in the observational processes, they affect each other in Lagrangian data assimilation. Therefore, the skill of the GB modes in the compressible flow case is lower than that in the incompressible flow case. This is also seen in Figure \ref{fig:covfig2}, which compares the posterior mean time series with the truth.

\begin{figure}[htbp]
	\centering
	\hspace*{-0.0cm}\includegraphics[width=1.00\textwidth]{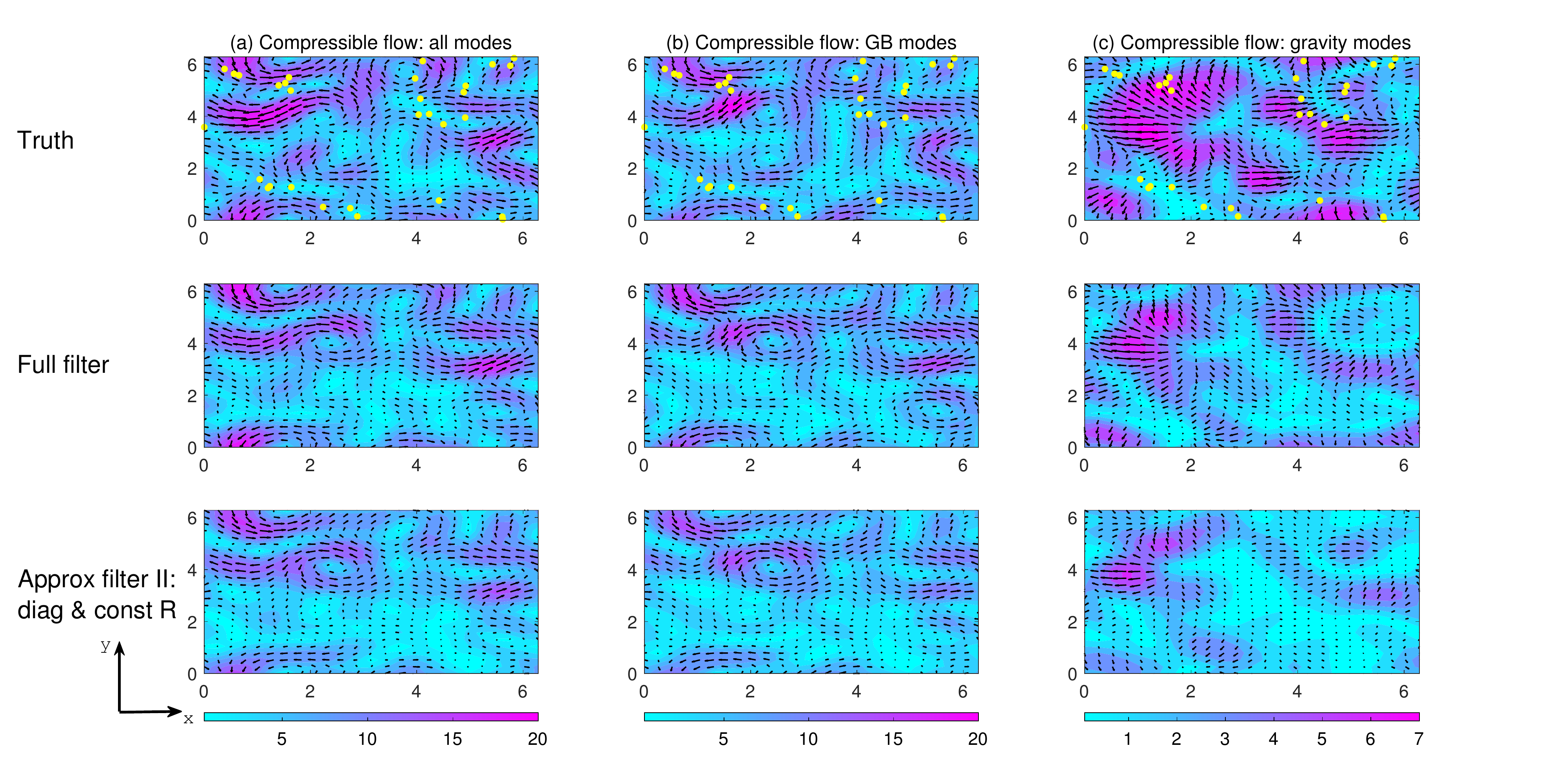}
	\caption{Comparison of the true flow field of the shallow water equation in physical space (top) and the recovered ones using the full filter (middle) and the Approximate filter II with a diagonal and constant covariance matrix (bottom). The yellow dots indicate the tracers’ locations, where $L=30$ tracers are utilized. The results correspond to $t=95$ (marked by a black square) in Figure \ref{fig:covfig2}, when the full filter recovers the truth quite well.}
	\label{fig:covfig3good}
\end{figure}

\begin{figure}[htbp]
	\centering
	\hspace*{-0.0cm}\includegraphics[width=1.00\textwidth]{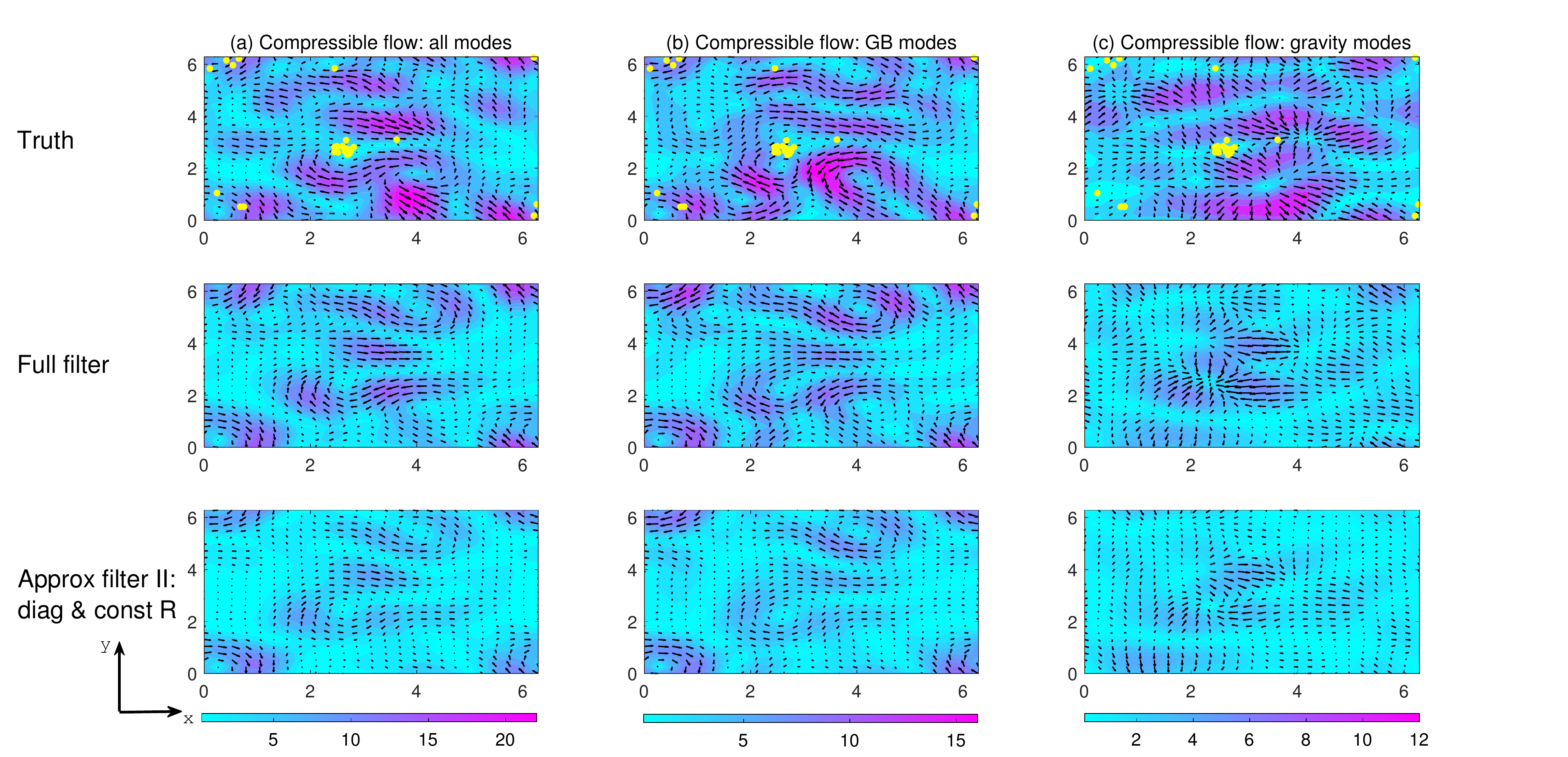}
	\caption{Similar to Figure \ref{fig:covfig3good} but the results correspond to $t=14$ (marked by a black square) in Figure \ref{fig:covfig2} when the full filter does not recover the truth accurately.}
	\label{fig:covfig3bad}
\end{figure}

To further understand the Lagrangian data assimilation skill in recovering the compressible flow field of the shallow water equation,  Fig.\ref{fig:covfig3good} and Fig.\ref{fig:covfig3bad} include two case studies. Both the full filter and the Approximate filter II recover the truth accurately at the time instant shown in  Fig.\ref{fig:covfig3good}. In contrast, the full filter becomes less accurate while Approximate filter II completely fails to capture the critical features in the case shown in  Fig.\ref{fig:covfig3bad}. The major source that causes the difference in these two cases is the strength of the compressible part of the flow. The ratio of the maximums between the GB and the gravity modes in Fig.\ref{fig:covfig3good} is 20:7 while that in Fig.\ref{fig:covfig3bad} is 15:12. When the gravity modes have a significant contribution to the overall flow field, the compressibility of the flow brings about the clustering of the tracers. As the tracers cannot cover the entire domain, the effective information provided by the tracers is reduced. This leads to the loss of skill in Lagrangian data assimilation utilizing the full filter. The Approximate filter II further suffers from such compressibility as it is derived based on the uniformly distributed tracer field resulting from incompressible flows.

\begin{figure}[htbp]
	\centering
	\hspace*{-0.0cm}\includegraphics[width=1.00\textwidth]{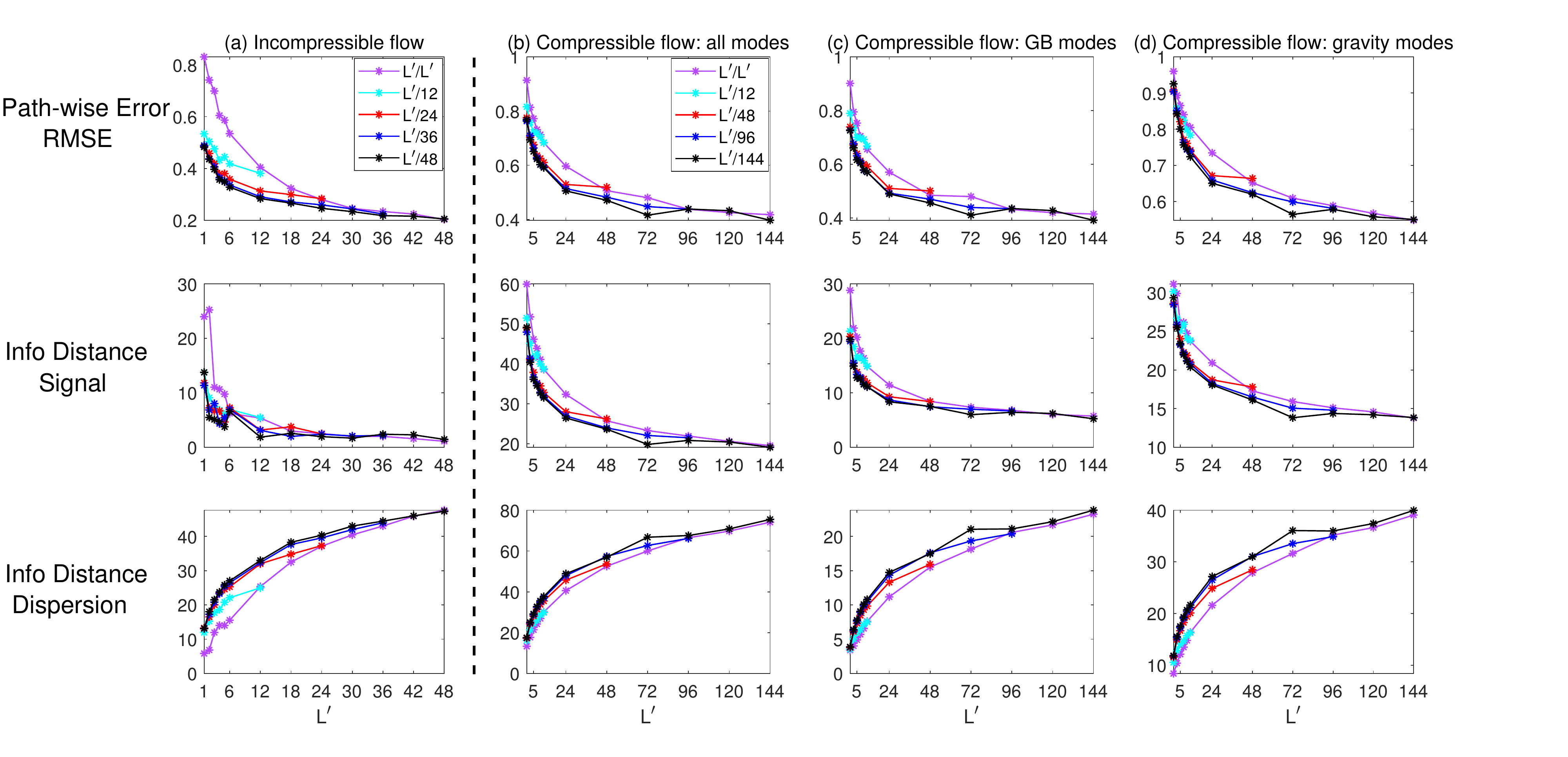}
	\caption{Skill scores using the Approximate filter III with a randomized selection of observations. In each panel, $L'/L$ means there are in total $L$ tracers in the domain, but only $L'$ randomly selected tracers are utilized to compute the posterior distribution at each instant. The selection of the $L$ tracers is taken at each time step. When $L'=L$, namely the last point in each curve, the approximate filter becomes the exact full filter (purple curve). Column (a) shows the case of the incompressible flow model, while Columns (b)--(d) show that of the compressible shallow water equation. The skill scores in Column (b) are computed based on all the shallow water equation modes, including both the GB and the gravity ones. The skill scores in Column (c) are computed only based on the filtered GB modes in the shallow water equation. Still, the Lagrangian data assimilation is applied to the entire shallow water equation. Column (d) is similar to Column (c) but for the filtered gravity modes. In this figure, the prefactor $\sqrt{L/L'}$ in \eqref{eq:prefactor} is not included for the illustration purpose since otherwise large error with $L'=1$ will dominate the error curve, and the details of the results cannot be seen. The case using the exact formula \eqref{eq:prefactor} for filtering is shown in Fig.\ref{fig:covfig4scale}. }
	\label{fig:covfig4}
\end{figure}

\begin{figure}[htbp]
	\centering
	\includegraphics[width=1\textwidth]{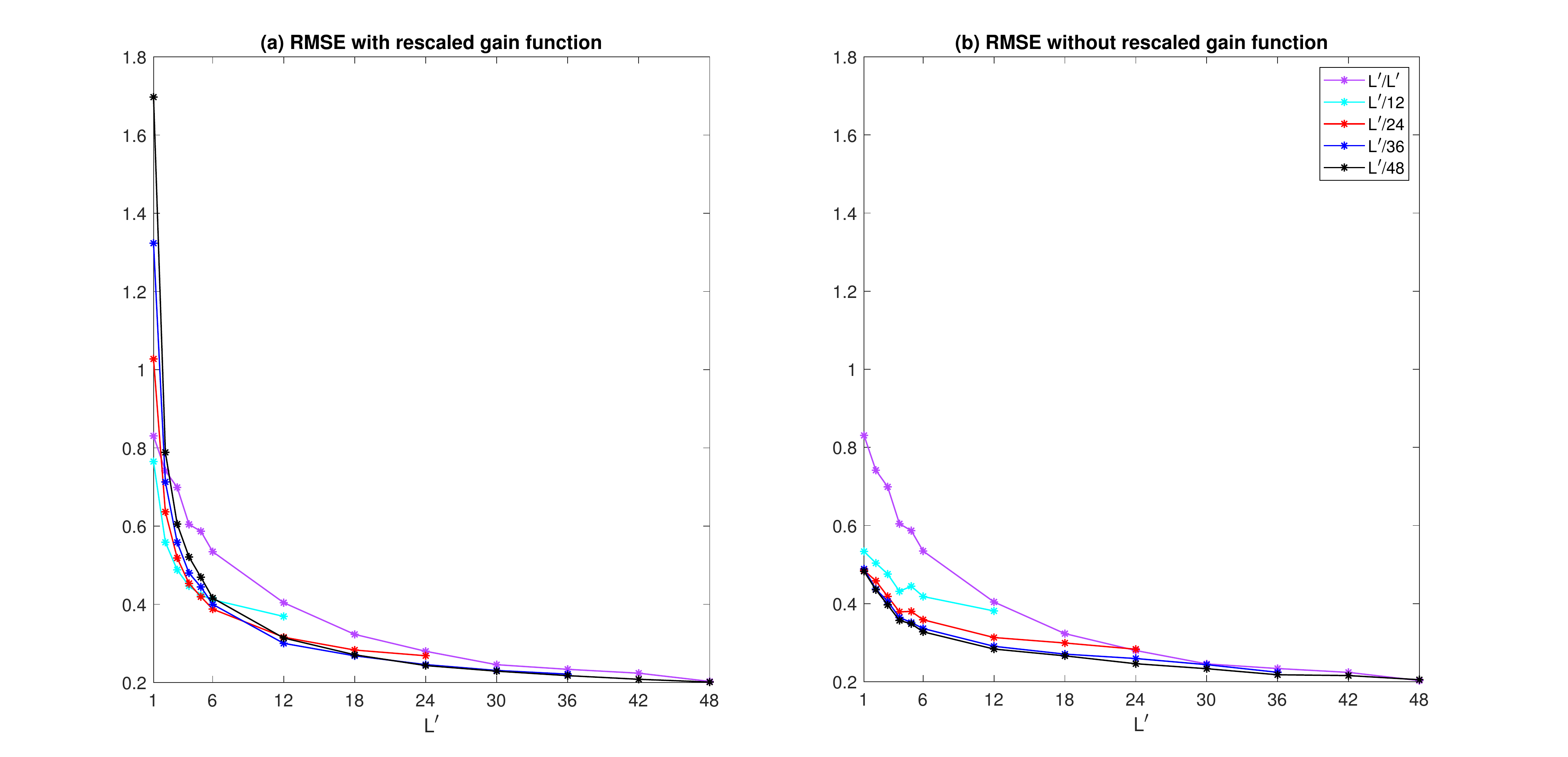}
	\caption{Comparison of the RMSE using randomly selected tracers with (Panel (a)) and without (Panel (b)) the prefactor in \eqref{eq:prefactor}. In each panel, $L'/L$ means there are in total $L$ tracers in the domain, but only $L'$ randomly selected tracers are utilized to compute the posterior distribution at each instant. The selection of the $L$ tracers is taken at each time step. When $L'=L$, namely the last point in each curve, the approximate filter becomes the exact full filter (purple curve), serving as a reference solution. The results here correspond to the incompressible flow field.}
	\label{fig:covfig4scale}
\end{figure}

%Finally,  Fig.\ref{fig:covfig4} (and Fig.\ref{fig:covfig4scale} in the Supplementary Material)
Finally,  Fig.\ref{fig:covfig4} and Fig.\ref{fig:covfig4scale}
present the skill scores using the Approximate filter III with a randomized selection of observations. It provides a numerical validation of Proposition \ref{theorem:random}, which says the number of tracers utilized in each instant can be significantly reduced if a randomized selection strategy is incorporated. Here the purple curve with legend $L'/L$ is served as a reference solution, which indicates the Lagrangian data assimilation skill if all the available $L$ tracers are used to recover the flow field. Therefore, the larger the gap between any curve with this one is, the more significant improvement the Approximate filter III brings. Note that, for the illustration purpose, the prefactor $\sqrt{L/L'}$ is removed in Fig.\ref{fig:covfig4} since otherwise, when $L'$ is very small with $L'=1$, the error is largely due to the failure of satisfying the assumption $L'\to\infty$. Nevertheless, once $L'\geq 3$, the Approximate filter III, both with or without such a prefactor, shows a significantly smaller error than the full filter using the same number of tracers. See Fig.\ref{fig:covfig4scale}.
When $L'$ is very small, the assumption $L'\to\infty$ in Proposition \ref{theorem:random} is broken, and therefore the filtering results are not expected to be accurate. In such a case, adding a big prefactor will amplify the random noise in the observational process. Thus, the filter without such a prefactor is even more skillful. Nevertheless, the advantage of such an approximate filter (both with and without the prefactor) becomes significant when $L'\geq 3$. Note that the total degree of freedom of the underlying incompressible flow field is 49, much larger than $L'$. This suggests the potential application of such an approximate filter in practice that requires using only a small number of randomly selected tracers to recover high-dimensional multiscale systems.
In addition, similar to the other two approximate filters, the method with the randomized selection of tracers works better for the incompressible case since the compressible flow suffers from the tracer cluster issue. Notably, the results indicate that such an approximate filter should also be helpful in many practical situations with discrete-in-time observations, provided that the observations arrive frequently.

\section{Comparison of Linear and Nonlinear Forecast Models in Lagrangian Data Assimilation}\label{Sec:ComparisonFilters}

As the underlying flow field is nonlinear in many practical problems, it is essential to understand the situations in which the efficient Lagrangian data assimilation algorithm using the approximate LSMs \eqref{eq:cgns}--\eqref{eq:filter} remains as skillful as the traditional ensemble-based data assimilation using the original nonlinear forecast model. To this end, the Lagrangian data assimilation framework \eqref{eq:cgns}--\eqref{eq:filter} is compared with the ensemble Kalman-Bucy filter (EnKBF) in recovering nonlinear turbulent flows. The EnKBF \cite{bergemann2012ensemble, de2018long} is an ensemble-based algorithm for continuous data assimilation applicable to nonlinear forecast models.

\subsection{A nonlinear flow model: the topographic barotropic model}
The barotropic quasi-geostrophic equations \cite{majda2006nonlinear} is an ideal model to study the complex nonlinear interaction between the zonal mean flow and the fluctuations via topography. The model reads:
\begin{equation}
\begin{aligned}
\frac{\partial q}{\partial t}+\nabla^{\perp}\psi\cdot \nabla q+u(t)\frac{\partial q}{\partial x}+\beta \frac{\partial\psi}{\partial x}&=f_q,\\
q=\Delta \psi+h,\qquad
\frac{du}{dt}-\fint h\frac{\partial \psi}{\partial x}&=-d_uu + f_u.
\end{aligned}
\end{equation}
Here the fluctuation is given in terms of the stream function $\psi$,
$\nabla^{\perp}\psi=(-\frac{\partial \psi}{\partial y},\frac{\partial \psi}{\partial x})$, $q$ is the potential vorticity, $u(t)$ is the zonal mean component (e.g., along $x$ direction), and the topography is given by $h = h(x,y)$. The parameter $\beta > 0$ is the contribution from the beta-plane effect, and $f_q$ and $f_u$ are external forcing terms. Here, $q$, $\psi$, and $h$ are all assumed to be $2\pi$-periodic functions in both variables $x$ and $y$ with zero average. The zonal mean velocity $u(t)$ is strongly coupled with the fluctuation, where the bar across the integral sign indicates that the integration has been normalized over the area of the domain.
Now consider the spectral formulation of the model with layered topography. The Fourier expansions of $\psi$ and $h$ are given by
\begin{equation}\label{eq:psi}
\psi(x,y,t)=\sum_{k\neq 0}\psi_k(t)e^{ik\mathbf{l}\cdot\mathbf{x}}\qquad\mbox{and}\qquad
h(x,y)=\sum_{k\neq 0}h_ke^{ik\mathbf{l}\cdot\mathbf{x}},
\end{equation}
where the expansion is along one characteristic wavenumber direction, which is denoted by $\mathbf{l}=(l_x,l_y)$ with $\lvert\mathbf{l}\rvert=1$. Therefore, the wavenumbers are given by $k\mathbf{l}$ with $k=\pm 1,\ldots, \pm K$. The spectral representation of the model is given by
\begin{equation}\label{eq:stoclayer}
\begin{aligned}
\frac{d \psi_k}{d t}=&-d_k\psi_k+ikl_x\bigg(\frac{\beta}{k^2|\mathbf{l}|^2}-u\bigg)\psi_k+i
\frac{kl_x}{k^2|\mathbf{l}|^2}h_ku+\sigma_k\dot{W_k},\\
\frac{d u}{d t}=&-d_uu-il_x\sum_{k\neq 0}kh_k\psi_k^*+\sigma_u\dot{W_k}
\end{aligned}
\end{equation}
where $\psi_k^*=\psi_{-k}$ and $h_k^*=h_{-k}$. Stochastic forcing is added to both the equations in \eqref{eq:stoclayer}. Without loss of generality $\mathbf{l} = (1, 0)$ is utilized here. Therefore, the velocity field in $y$ direction is completely determined by $\psi$ with $v = \partial\psi/\partial x$ while the velocity in $x$ direction is  $u$. The topographic effect is given by $h_1=H_1/2-iH_1/2$, $h_2=H_2/2-iH_2/2$, and $h_k=\frac{\sin\theta_k}{4k^p}-i\frac{\cos\theta_k}{4k^p}$ for $k\geq 3$, where $\theta_k$ is a random phase and $p$ is a power that controls the effects of the small-scale topography. The following parameters are utilized in the experiments $d_k=d_u=0.0125$, $H_1=1$ and $H_2=1/2$. Two regimes are considered here:
\begin{itemize}
	\item \emph{Regime I: Nearly Gaussian regime:} The fluctuation modes are subject to relatively strong noise forcing with strength  $\sigma_{v,k}=1/\sqrt{2}$,   $\sigma_{v,k}=1/(2\sqrt{2})$, $\sigma_{v,k}=1/(\sqrt{2}k^2)$ for $k=3,\ldots,K$ while the noise strength in the zonal mean flow is $\sigma_{u}=1/2\sqrt{2}$.
	\item \emph{Regime II: Strongly non-Gaussian regime:} The zonal mean flow is strongly forced with white noise strength $\sigma_{u}=1/\sqrt{2}$ while relatively weak noises $\sigma_{v,1}=\sigma_{v,2}=1/4\sqrt{2}$ and $\sigma_{v,k}=1/(2\sqrt{2}k^2)$ for $k=3,\ldots,K$ are added to the fluctuation modes.
\end{itemize}
For the tests in both regimes, $K=6$ is utilized. The standard Runge-Kutta 4th order method is utilized for the numerical integration of the deterministic part, while the numerical integration of the stochastic part is the same as the Euler-Maruyama scheme. The numerical integration time step is $\Delta{t} = 0.001$. A long realization with $400$ time units is utilized as the reference solution to compute the skill score in Lagrangian data assimilation.

\begin{figure}[htbp]
	\centering
	\includegraphics[width=1\textwidth]{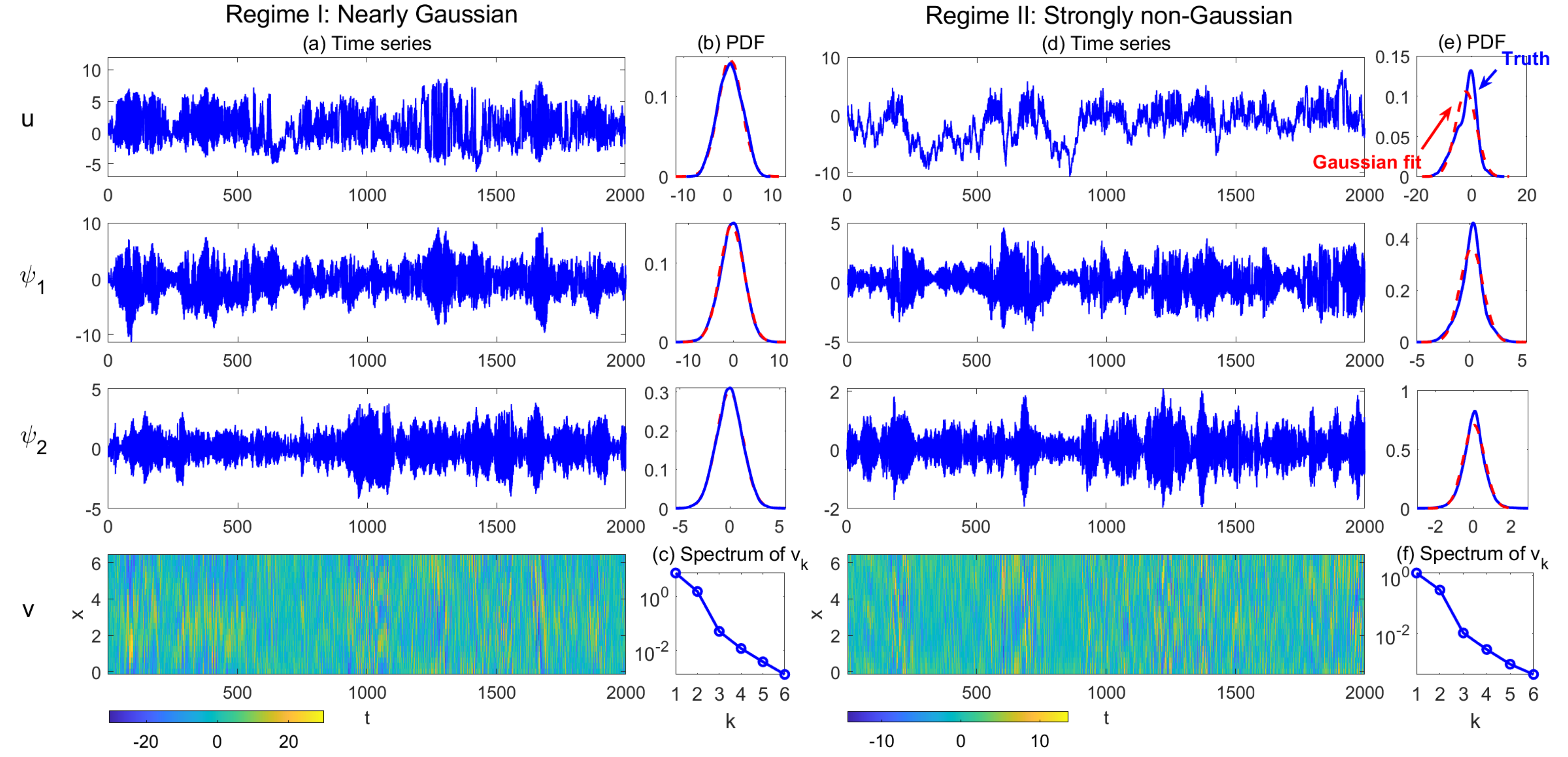}
	\caption{Dynamical regimes of the layered topographic model \eqref{eq:stoclayer}. Panels (a)--(c) are for the nearly Gaussian regime. Panel (a): time series of $u$, $\psi_1$ and $\psi_2$ as well as the spatiotemporal reconstruction of the velocity field $v$. Panel (b): the PDFs of  $u$, $\psi_1$ and $\psi_2$ (blue) together with their Gaussian fit (red). Panel (c): the energy spectrum of $v_k$. Panel (d)--(f) are similar to Panels (a)--(c) but are for the non-Gaussian regime. }
	\label{fig:layerfig1}
\end{figure}

Fig.\ref{fig:layerfig1} shows the dynamical and statistical features of these two regimes. In the non-Gaussian regime, the large white noise forcing in $u$ excites a strong competition between two alternating states: a highly intermittent flow field $v$ when $u>0$ and a nearly steady flow when $u<0$. The competition between two alternating states leads to intermittent eastward zonal flow $u>0$ and steady westward flow $u<0$. The intermittent nature of the flow field triggers non-Gaussian fat-tailed distributions of $v$. On the other hand, in the nearly Gaussian regime, strongly multiscale features emerge in the time series of all the variables, which include multiple fast scales with rapid oscillations and a slowly varying long-term tendency. In particular, there are two distinct dominant frequencies in the fast oscillations. The extremely fast oscillation appears when the zonal flow goes steadily with $u<0$, while the moderately fast one occurs when the zonal flow becomes intermittent with an average velocity $u>0$. In contrast to the non-Gaussian regime, the nearly Gaussian statistics in this regime are due to the comparable amplitude of $v_k$  at different time instants. %{ Fig.\ref{fig:layerfig1_zoomin} and Fig.\ref{fig:layerfig2_zoomin} show the enlarged versions of Fig.\ref{fig:layerfig1} that include the times series of the large-scale zonal velocity $u$, the real part of the first and the second Fourier modes $\psi_1$ and $\psi_2$, and the entire fields of $v$ and $T$ in physical space as a function of time. With the zoom-in details, it can be seen that, in the non-Gaussian regime, the competition between two alternating states leads to intermittent eastward zonal flow $u>0$ and steady westward flow $u<0$. On the other hand, in the near-Gaussian regime, strongly multiscale structures appear in the time series of all variables together with highly oscillating fast scales and a slow varying long-time structure.}
% See also the enlarged figures in the Supplementary Material.

%\begin{figure}[htbp]
%	\centering
%	\includegraphics[width=1\textwidth]{Figs_all/SI_Topography_1.pdf}
%	\caption{Time series of $u$, $\psi_1$ and $\psi_2$ as well as the spatiotemporal reconstruction of the velocity field $v$. }
%	\label{fig:layerfig1_zoomin}
%\end{figure}

%\begin{figure}[htbp]
%	\centering
%	\includegraphics[width=1\textwidth]{Figs_all/SI_Topography_2.pdf}
%	\caption{Time series of $u$, $\psi_1$ and $\psi_2$ as well as the spatiotemporal reconstruction of the velocity field $v$. }
%	\label{fig:layerfig2_zoomin}
%\end{figure}
\subsection{Comparison of the Lagrangian data assimilation skill}
The layered topographic model \eqref{eq:stoclayer} is nonlinear. Therefore, the Lagrangian data assimilation using the LSMs contains model errors. Knowing the exact nonlinear model here facilitates the quantification of model error and uncertainty in affecting the data assimilation skill. In addition, according to \eqref{eq:stoclayer}, the interaction between different $\psi_k$ is only through the zonal flow $u$. This allows a straightforward way to build an approximate filter with model reduction that nevertheless preserves the model structure.
The following $5$ filters are utilized here:
\begin{enumerate}
	\item \emph{EnKBF with all modes.} This is treated as the reference solution, where the forecast model in the Lagrangian data assimilation is the perfect model \eqref{eq:stoclayer}.
	\item \emph{EnKBF with $u+v_{\pm 1}+v_{\pm 2}$.} This filter includes  model error due to a dimension reduction. Only the equations of $u$ and the those of $\psi_{\pm1}$ and $\psi_{\pm2}$ (equivalently to $v_{\pm1}$ and $v_{\pm2}$) are utilized as the forecast model. The goal of such a filter is to efficiently recover the large-scale features of the flow field.
	\item \emph{LSMs with all modes.} The original layered topographic model \eqref{eq:stoclayer} is approximated by the LSMs in Lagrangian data assimilation using the framework \eqref{eq:cgns}.
	\item \emph{LSMs with $u+v_{\pm 1}+v_{\pm 2}$.} Only the true signals of $u$, $\psi_{\pm1}$ and $\psi_{\pm2}$ from the original layered topographic model \eqref{eq:stoclayer} are utilized to calibrate the LSMs as the forecast model in Lagrangian data assimilation  \eqref{eq:cgns}. This filter contains model error from both dimension reduction and linear approximation of the forecast model.
	\item \emph{Purely data-driven LSMs.} Here, the ansatz of the perfect model is assumed to be completely unknown. Therefore, the LSMs are developed purely based on the given spatiotemporal data of the two-dimensional velocity field. The LSMs here contain two wavenumbers $k_x$ and $k_y$. Both range from $-K$ to $K$.
\end{enumerate}
Algorithm \ref{Alg:ParameterEstimation} is adopted to estimate the parameters in all the LSMs mentioned above.

\begin{figure}[htbp]
	\centering
	\includegraphics[width=1\textwidth]{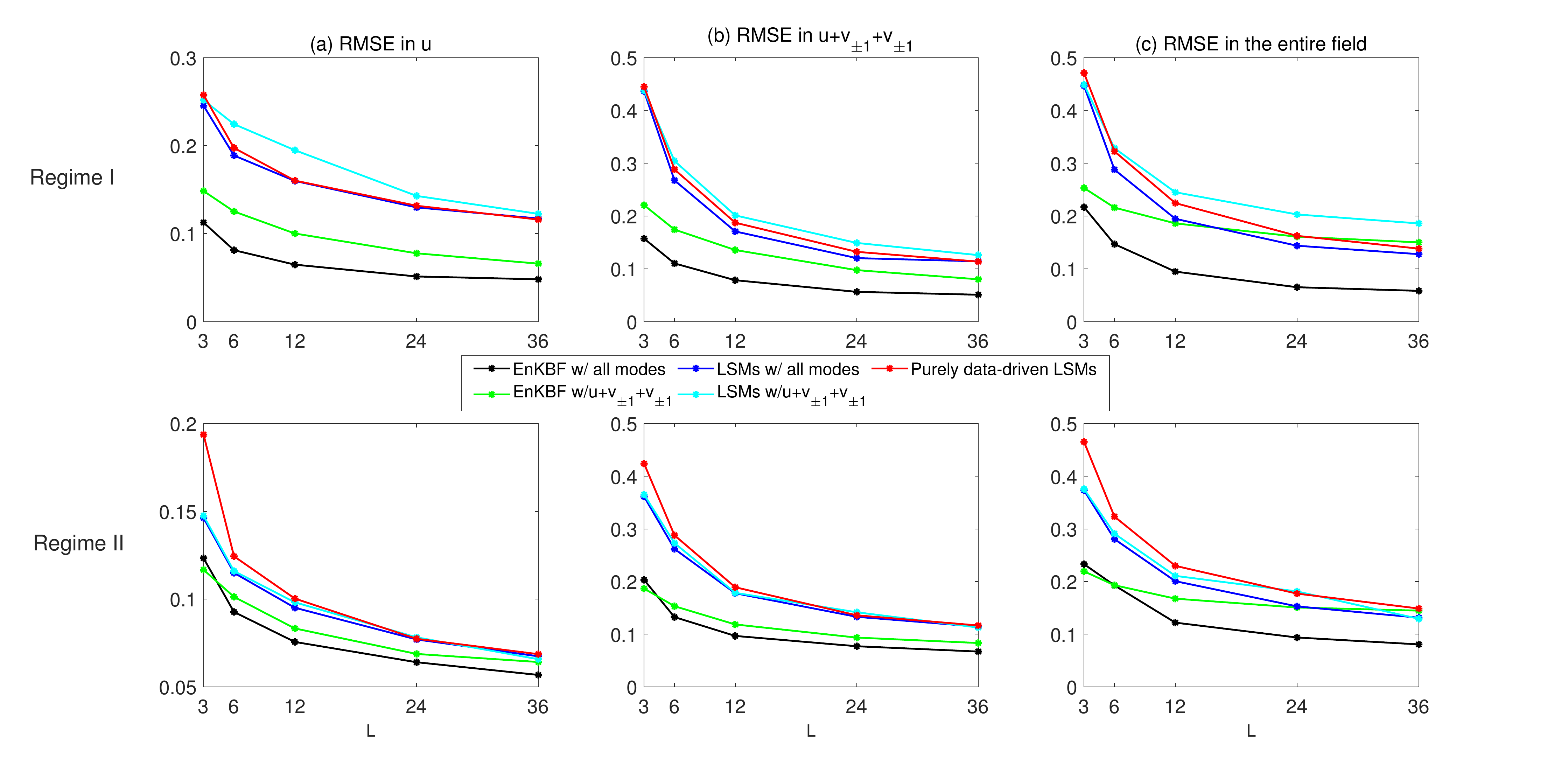}
	\caption{RMSE in the posterior mean estimate in filtering the layered topographic model. The black curve shows the EnKBF using the full model. The green curve shows the EnKBF using the approximate model that includes only $u$ and $v_{\pm1}$ and $v_{\pm2}$. The blue curve shows the data assimilation using the LSMs to approximate all the modes. The cyan curve shows the data assimilation using the LSMs but only includes $u$ and $v_{\pm1}$ and $v_{\pm2}$. The red curve shows the data assimilation using the LSMs, where the ansatz of the perfect model is assumed to be unknown, and therefore the LSMs contain two wavenumbers $k_x$ and $k_y$. Both range from $-K$ to $K$.}
	\label{fig:layerfig2}
\end{figure}

Fig.\ref{fig:layerfig2} shows the RMSE in the posterior mean estimate as a function of $L$. The EnKBF with all modes (black) is always the best among different filters since the perfect model is utilized. The EnKBF with $u+v_{\pm 1}+v_{\pm 2}$ (green) is also accurate in recovering the leading a few modes. Although the three filters with LSMs are not as skillful as the EnKBFs, the RMSE is below $30\%$ for all the modes when more than $6$ tracers are utilized in both regimes, which justifies the use of such approximate filters. In addition, the relatively large error, when only a small number of tracers is utilized, is more attributed to the uncertainty in the parameter estimation.
It is worth highlighting that the filters with LSMs do not illustrate a significant deterioration in the non-Gaussian regime. Although the marginal distribution (e.g., the flow field distribution) is strongly non-Gaussian, the conditional distribution of the flow given the observations may not show a strong non-Gaussianity at all time instants. This is seen from the time series in Panel (d) of Fig.\ref{fig:layerfig1}, where the non-Gaussian distribution of $\psi_1$ (and $\psi_2$) is made by a Gaussian mixture with two components with different variances. Thus, at each instant, the distribution of $\psi_1$ (and $\phi_2$) behaves more towards Gaussian. A similar feature can be observed in $u$. The forecast distribution of $u$ within the short term does not seem to have a bimodal behavior. These findings allow the framework \eqref{eq:cgns} to have the potential to be applied to many other non-Gaussian flows.

\begin{figure}[htbp]
	\centering
	\includegraphics[width=1\textwidth]{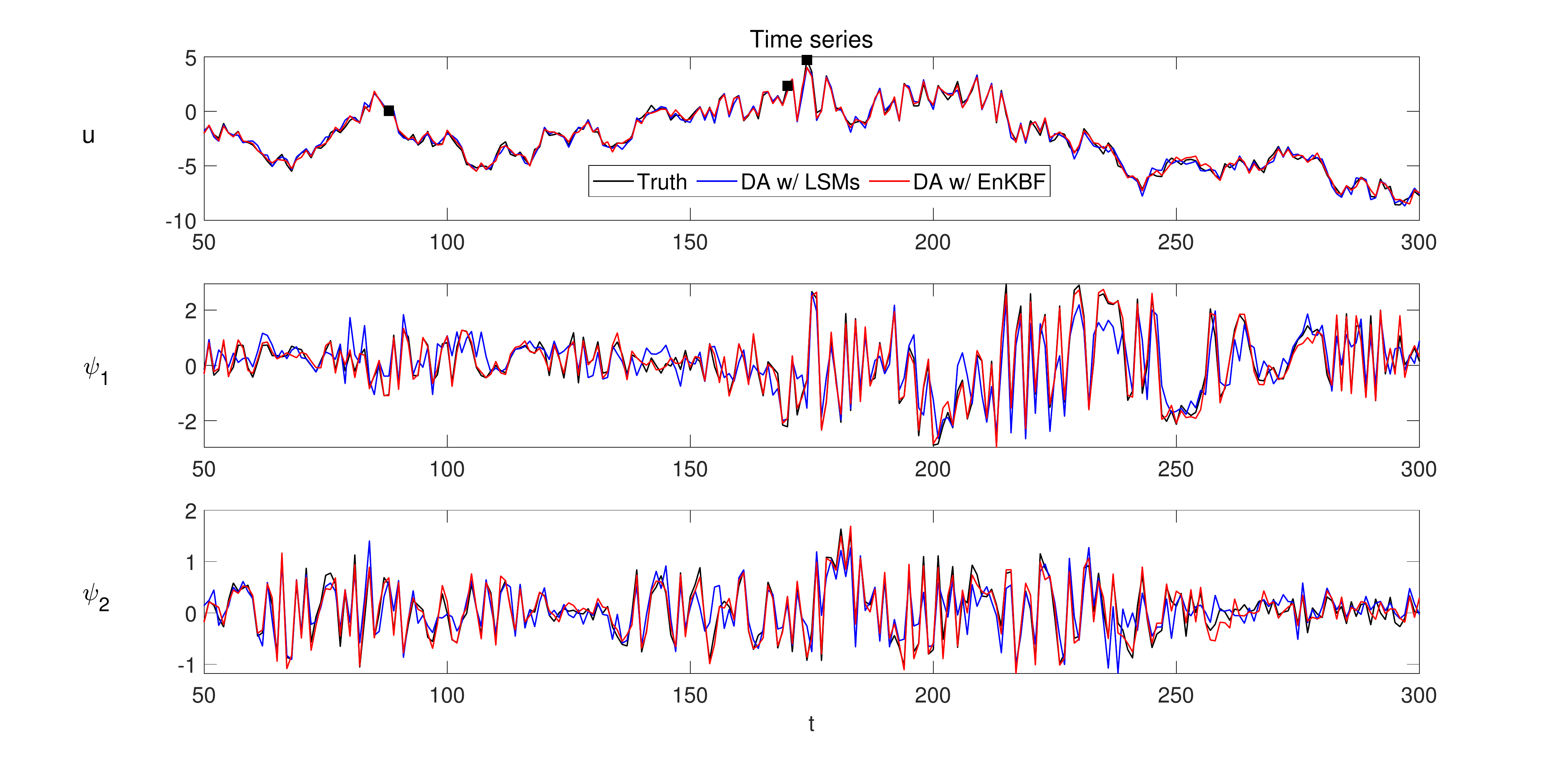}
	\caption{Comparison of the truth with the posterior mean estimate using the filter with the LSMs and the EnKBF in the non-Gaussian regime. Both filters include all the modes, $L=6$.}
	\label{fig:layerfig3}
\end{figure}

\begin{figure}[htbp]
	\centering
	\includegraphics[width=1\textwidth]{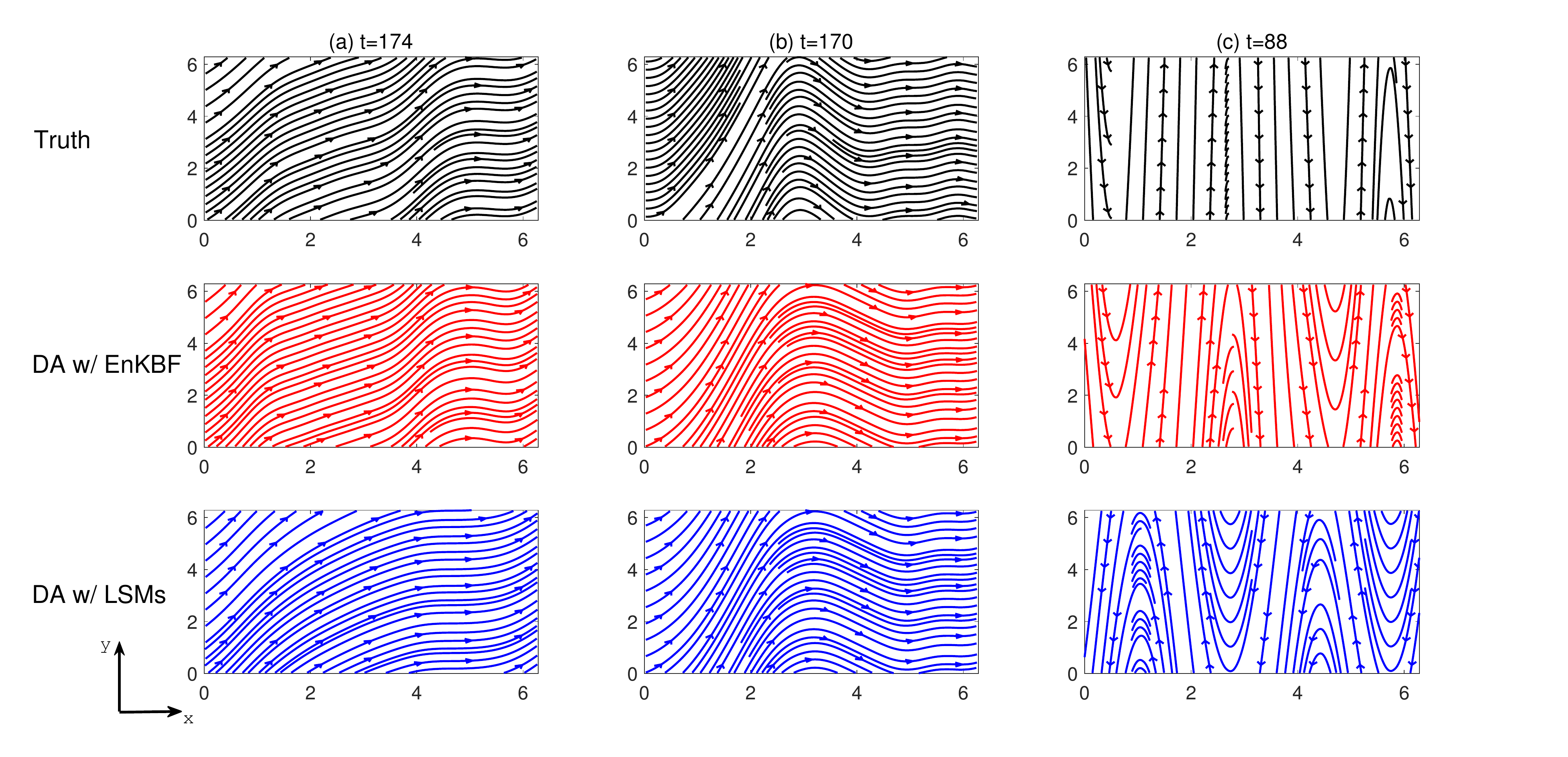}
	\caption{The reconstructed spatial pattern at three different time instants marked by the black squares in Fig.\ref{fig:layerfig3}, $L=6$.}
	\label{fig:layerfig4}
\end{figure}

Fig.\ref{fig:layerfig3} compares the time series of the truth with the posterior mean estimate using the filter with the LSMs and the EnKBF containing all the modes. It is seen that the filter with the LSMs recovers the true signal quite accurately and is only slightly worse than the EnKBF. Fig.\ref{fig:layerfig4} shows the reconstructed spatial pattern at three different time instants marked by the black squares in Fig.\ref{fig:layerfig3} with $L=6$. The recovered entire velocity fields at both $t=174$ and $t=170$, where the zonal flow $u$ is strong or moderate, are very accurate. On the other hand, when $u$ is near zero, then the overall amplitude of the velocity field is much weaker and is dominated by the velocity in the $y$ direction. In such a case, errors in the recovered field are observed, especially in the filter with the LSMs.

\subsection{The Purely Data-Driven LSMs}
\begin{figure}[htb!]
	\centering
	\includegraphics[trim={0cm .0cm 0cm 0cm},clip,width=6.0in]{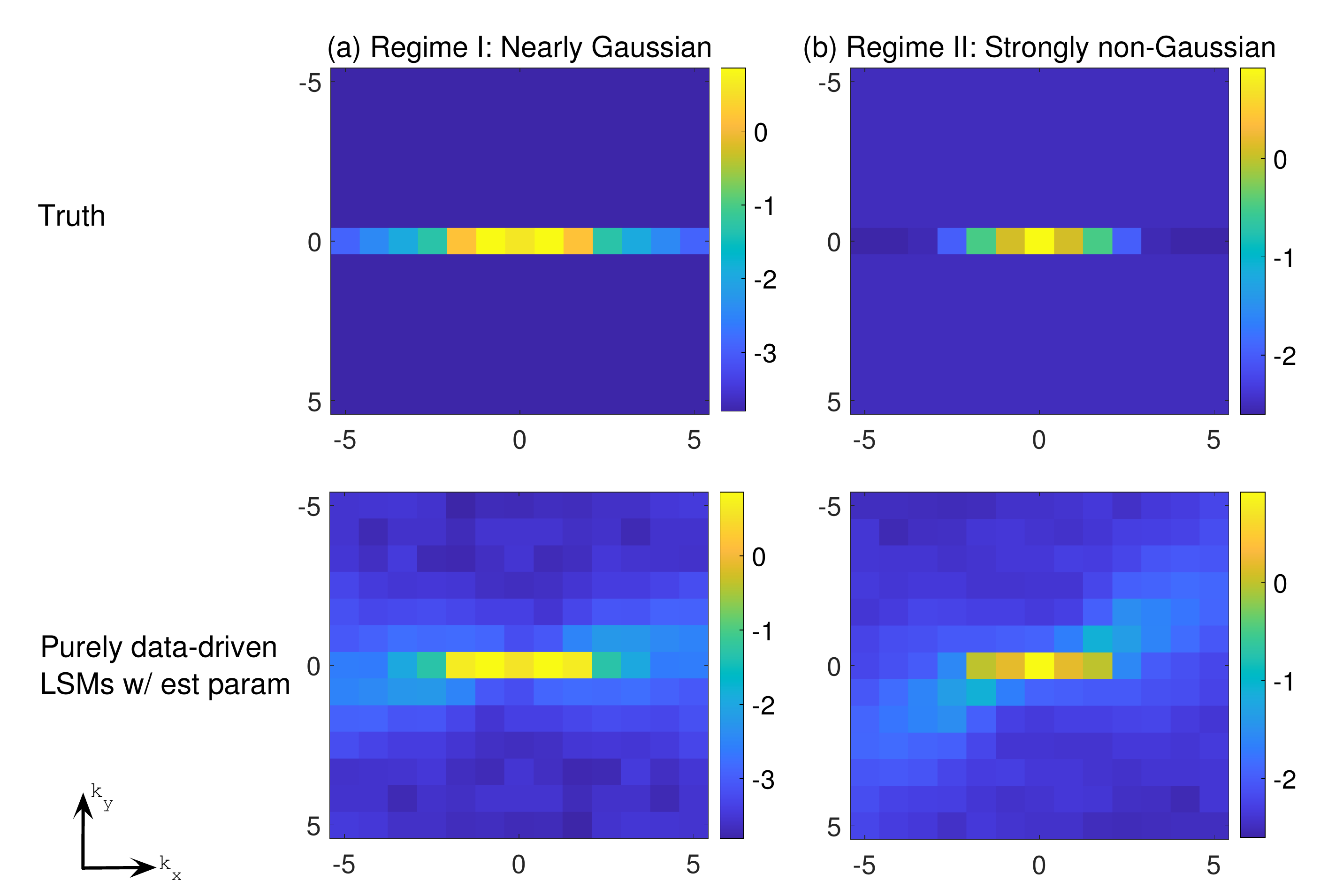}
	\caption{Comparison of the true energy spectrum and the  recovered energy spectrum from the Lagrangian data assimilation with the data-driven LSMs, where $K=6$ and $L=6$ tracers are utilized for the parameter estimation and the Lagrangian data assimilation. The energy is shown in the log scale to better display the small values.  }
	\label{fig:layfigsi1}
\end{figure}

This section provides additional results of using the purely data-driven LSMs to perform the Lagrangian data assimilation for the
layered topographic model. In the data-driven LSMs, the ansatz of the model that generates the trajectories is unknown. The only information available is the trajectories of the tracers. Therefore, even though in the layered topographic model all the wavenumbers are along a specific direction, the LSMs utilized here contain two wavenumbers $k_x$ and $k_y$, ranging from $-K$ to $K$. In the following, two cases are considered, $K=6$ and $K=7$. The former case is referred to as the standard data-driven LSMs, as it has the same truncation of the wavenumbers up to $K=6$. The latter does not even utilize the information of the total wavenumbers in generating the true signal. In both cases, the parameters in the LSMs are estimated from the observed Lagrangian tracers, as was described in Algorithm \ref{Alg:ParameterEstimation}. Fig.\ref{fig:layfigsi1} compares the true and recovered energy spectrum in log scale with the data-driven LSMs. Note that the logarithm of the energy is taken to display the spectrum better since, otherwise, the small values are hard to be seen from the figure. The figure indicates that the energy recovered from the Lagrangian data assimilation concentrated on the trace $k_y=0$, which illustrates the skill of the data-driven LSMs in parameter estimation and data assimilation.

\begin{figure}[htb!]
	\centering
	\includegraphics[trim={0cm .0cm 0cm 0cm},clip,width=6.5in]{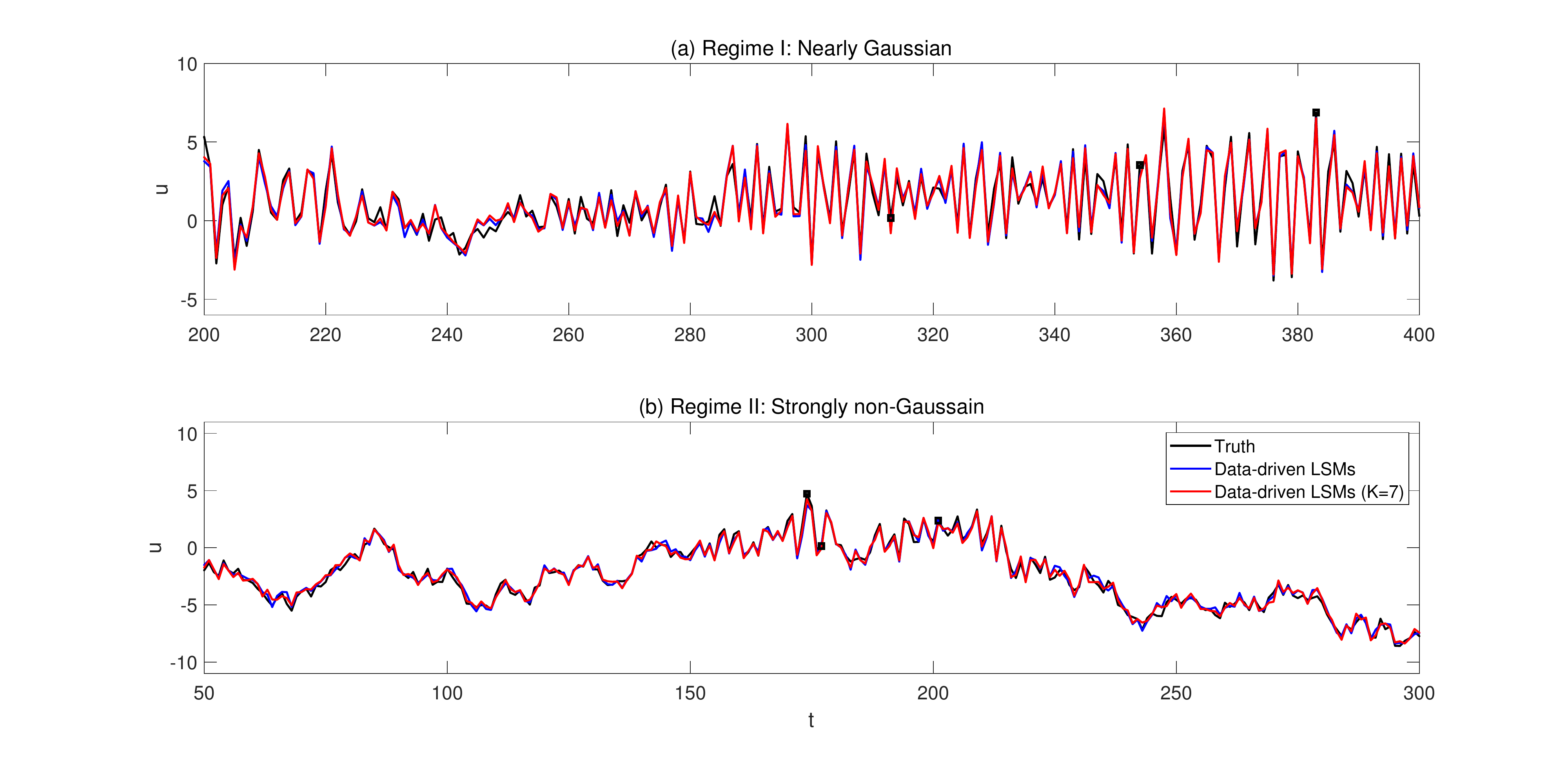}
	\caption{ Comparison of the truth with the posterior mean estimate using the standard data-driven LSMs  and the data-driven LSMs
		with $K=7$ and $L=6$. The three black squares correspond to the time instants of showing the spatial patterns in Fig.\ref{fig:layfigsi2_regI} and Fig.\ref{fig:layfigsi2_regII}}
	\label{fig:layfigsi4}
\end{figure}

Fig.\ref{fig:layfigsi4} compares the time series and the recovered $u$ with the data-driven LSMs, with $L=6$ tracers being utilized.  The data-driven LSMs can accurately recover the zonal mean flow $u$  even with a small number of tracers and little prior information about the spectrum domain.

\begin{figure}[htb!]
	\centering
	\includegraphics[trim={0cm .0cm 0cm 0cm},clip,width=6.5in]{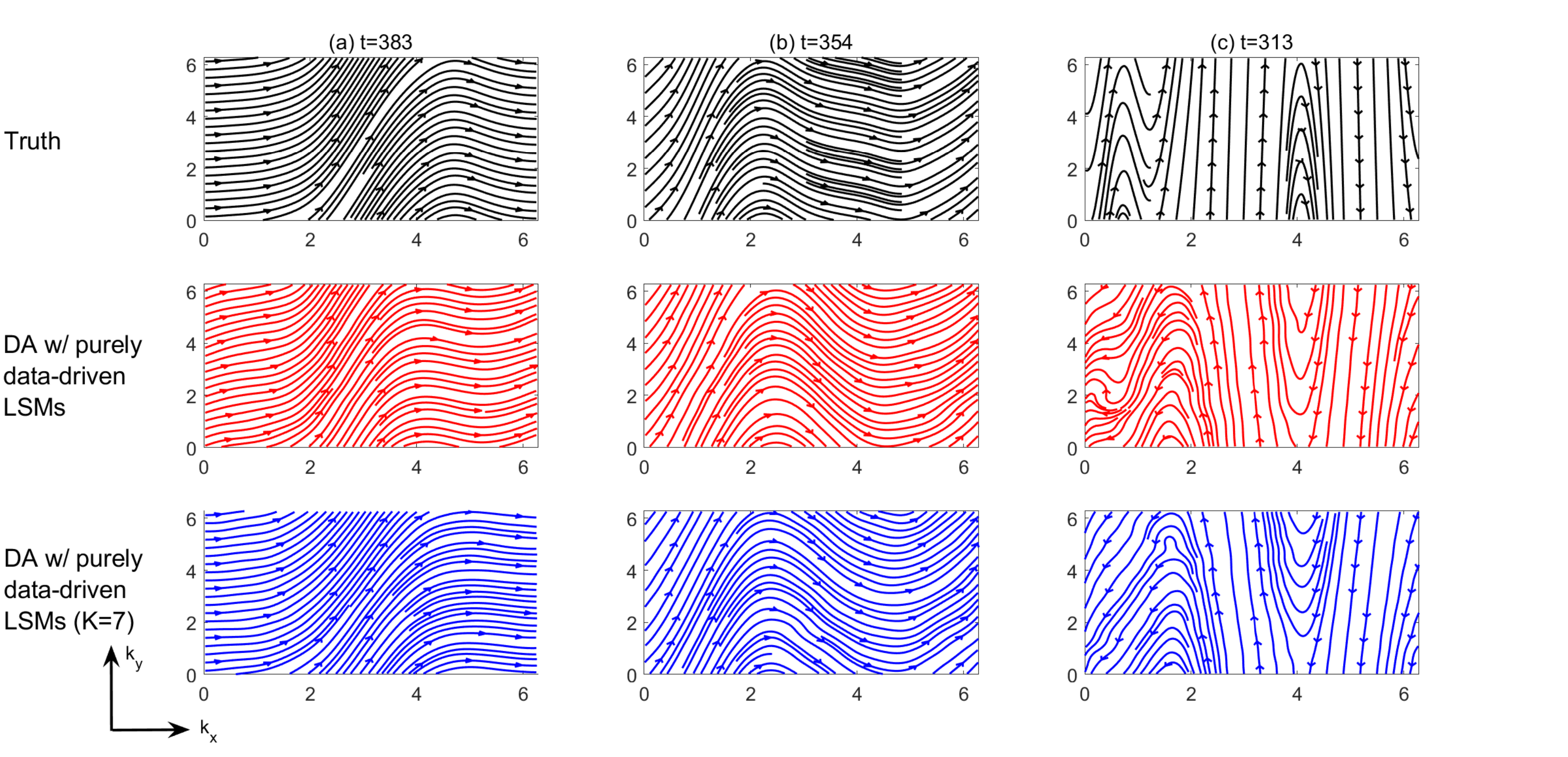}
	\caption{The reconstructed spatial pattern of the nearly Gaussian regime at three different time instants marked by the black squares in
		Fig.\ref{fig:layfigsi4}, where  $L=6$.}
	\label{fig:layfigsi2_regI}
\end{figure}

\begin{figure}[htb!]
	\centering
	\includegraphics[trim={0cm .0cm 0cm 0cm},clip,width=6.5in]{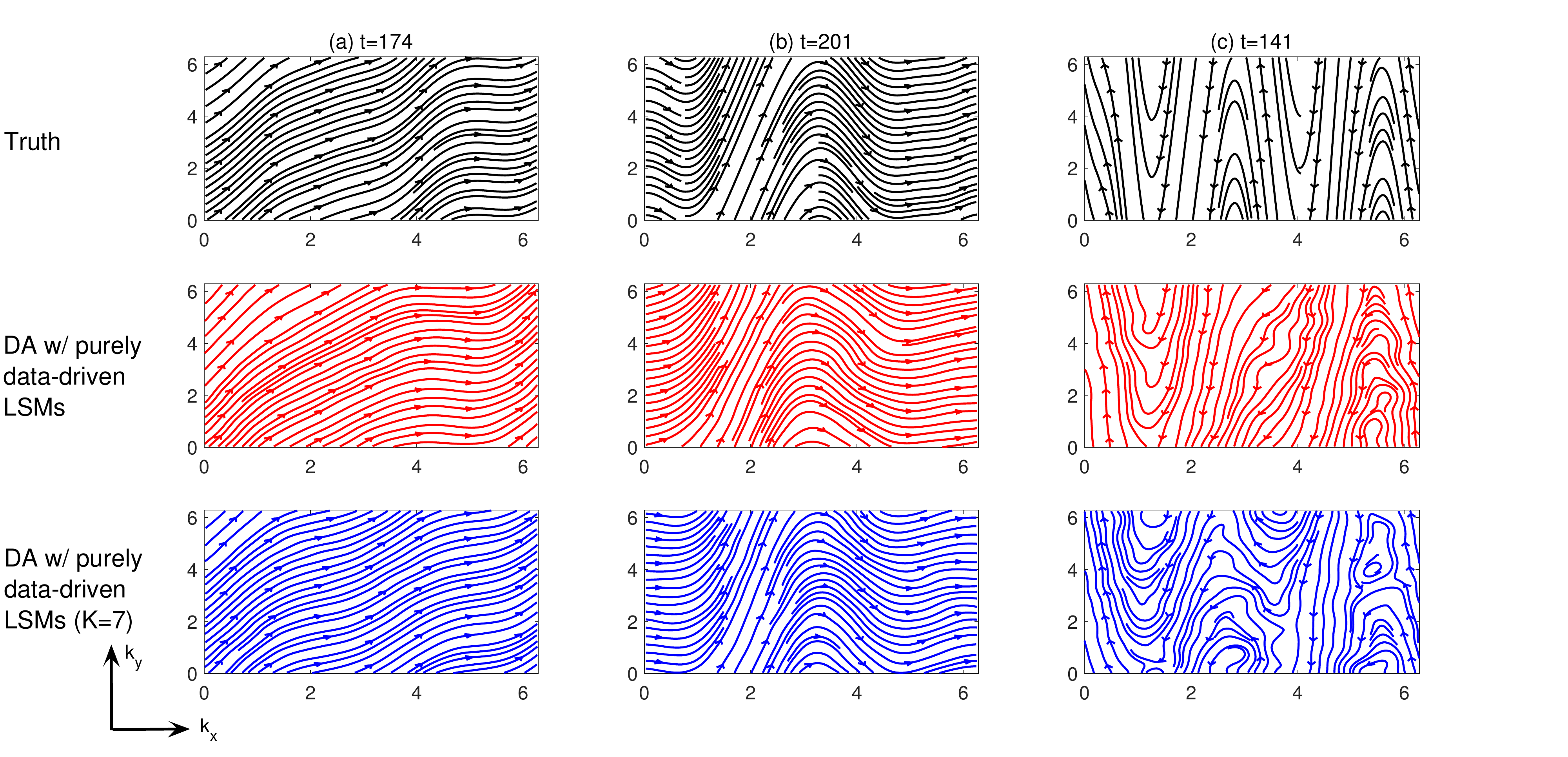}
	\caption{ The reconstructed spatial pattern of the strongly non-Gaussian regime at three different time instants marked by the black squares in
		Fig.\ref{fig:layfigsi4}, where  $L=6$.}
	\label{fig:layfigsi2_regII}
\end{figure}

\begin{figure}[htb!]
	\centering
	\includegraphics[trim={0cm .0cm 0cm 0cm},clip,width=5.0in]{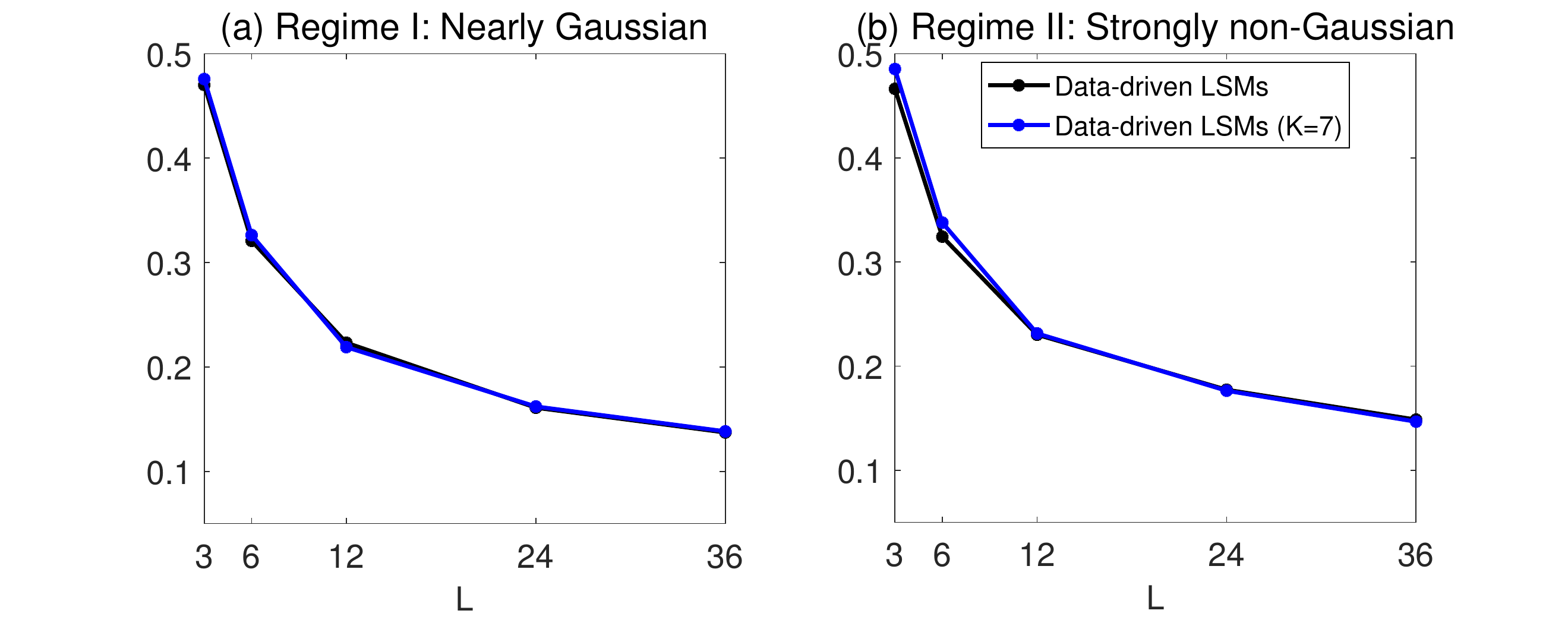}
	\caption{The RMSE in the posterior mean estimate in filtering the layered topographic model with the data-driven LSMs.}
	\label{fig:layfigsi3}
\end{figure}

Fig.\ref{fig:layfigsi2_regI} and Fig.\ref{fig:layfigsi2_regII} compare the reconstructed spatial pattern with the data-driven LSMs at three different time instants marked by the black squares in Fig.\ref{fig:layfigsi4}. For the nearly-Gaussian regime, the recovered entire velocity fields show high accuracy at both $t = 383$ and $t = 354$, where the zonal flow $u$ is strong or moderate. On the other hand, like other filters, when $u$ is near zero, the velocity field is dominated by the component in the $y$ direction, and the overall amplitude of the velocity field is weak. As a consequence, some errors are observed. Nevertheless, the overall error remains small even with no prior information about the spectrum space, e.g., assuming $K=7$ in the LSMs. It is also worthwhile to mention that the non-Gaussian regime behaves similarly to the nearly Gaussian regime, which again indicates the skill of the LSMs in recovering the non-Gaussian features with the help of the Lagrangian tracers. Fig.\ref{fig:layfigsi3} displays the RMSE in the posterior mean estimate in filtering the layered topographic model with the data-driven LSMs, which confirms the above conclusions.

\section{Conclusion}\label{Sec:Conclusion}
An analytically tractable mathematical framework for continuous-in-time Lagrangian data assimilation is developed to study the skill of nonlinear Lagrangian data assimilation using linear stochastic forecast models. The framework advances the development of an efficient parameter estimation algorithm (Section \ref{Sec:ParameterEstimation}) and the uncertainty quantification of several approximate filters (Section \ref{Sec:ReducedFilters}). The diagonal approximation of the posterior covariance is shown to be a suitable strategy for reducing the computational cost of the filters in dealing with incompressible flows, which appear in many ocean science applications. On the other hand, randomly selecting a small number of tracers at each time step as observations can reduce the computational cost while retaining the data assimilation accuracy. Such a finding, in particular, provides crucial justifications and guidelines for adopting non-interacting sea ice floe trajectories as Lagrangian tracers to recover the underlying ocean flow field, where these floe trajectories often appear randomly in the Arctic region. They typically last for a shorter time than the ocean's advective time scale before contracting with others \cite{covington2022bridging}.
The framework also allows a systematic study of model error for filtering nonlinear turbulent flows using the linear stochastic forecast models (Section \ref{Sec:ComparisonFilters}). It is shown that the linear stochastic forecast models are skillful even for filtering the layered topographic flow field with non-Gaussian features.
{{It is worthwhile to remark that although all the test examples considered in this work have periodic boundary conditions, the general framework developed here is not limited to such a boundary condition. The only major modification will be that $g_\mathbf{k}(\mathbf{x})$ will no longer be a Fourier basis function. Nevertheless, the data assimilation approach developed here can still be carried out as it is applied to the time series of the coefficients associated with the basis functions.}
Future work includes applying the framework to more complicated systems and incorporating it with other reduced-order modeling or filtering strategies.}

\section*{Acknowledgments}
N.C. is partially funded by ONR N00014-21-1-2904 and N00014-19-1-2421. S.F. was a former postdoc research associate under the second grant. The authors sincerely thank Dr. Andrew Stuart for providing many useful and crucial suggestions.
\section*{Data Availability Statement }
{
The Matlab scripts used in this work are posted on Githb repository (https://github.com/aggietx/source-code-for-uq-with-lsm).
}

\section*{Appendix}
\subsection{The linear shallow water system}
This section includes the details of the shallow water model utilized in Section \ref{Sec:ReducedFilters} of the main text for generating the compressible flow field. The linear shallow water equation includes incompressible geophysically balanced (GB) and compressible gravity modes.
It starts from a linear rotating shallow water equation,
\begin{equation}\label{SI_SWE}
\begin{aligned}
\qquad\qquad\frac{\partial \mathbf{u}}{\partial t}+\mbox{Ro}^{-1}\mathbf{u}^\bot&=-\mbox{Ro}^{-1}\nabla \eta,\\
\qquad\qquad\frac{\partial \eta}{\partial t}+\mbox{Ro}^{-1}\delta\nabla\cdot \mathbf{u}&=0,
\end{aligned}
\end{equation}
where $\mathbf{u}$ is the two-dimensional velocity field and $\eta$ is the height function in the domain $[0,2\pi)^2$ with double periodic boundary conditions.
In \eqref{SI_SWE}, $\mbox{Ro}$ is the Rossby number, a non-dimensional number representing the ratio between the Coriolis term and the advection term. On the other hand, Define Fr as the Froude number and therefore $\delta=\mbox{Ro}^2\mbox{Fr}^{-2}$. As in many applications, $\mbox{Fr}$ is set to be the same value as $\mbox{Ro}$ here.
%  For most atmosphere and ocean problems, $\epsilon$ ranges from $O(10^{-1})$ to $O(1)$, representing fast to moderate rotation while    $\delta=1$ is a typically choice.
The velocity field ${\bf u}$ and the geophysical height function $\eta$ can be written as
\begin{equation}\label{SI_SWE_Velocity_field}
\left(
\begin{array}{c}
\mathbf{u}(\mathbf{x}, t) \\
{\eta}(\mathbf{x}, t) \\
\end{array}
\right)
= \sum_{\mathbf{k} \in \mathbb{Z}^2, \alpha\in \{B,\pm\}}\widehat{v}_{\mathbf{k},\alpha}(t) e^{i \mathbf{k} \cdot \mathbf{x}} \mathbf{r}_{\mathbf{k},\alpha},
\end{equation}
% where $\mathbf{u}_o$ is the two-dimensional velocity field and $\eta$ is the height function.
The wavenumber $\mathbf{k}=(k_x,k_y)$. The set $\{B,+,-\}$, where $\alpha$ belongs, contains three elements.

In \eqref{SI_SWE_Velocity_field}, the modes with $\alpha = B$ represent the incompressible GB modes,  in which the GB relation $\mathbf{u}^\perp = -\nabla \eta$ always holds. The GB flows are non-divergent, which is embodied in the eigenvector $\mathbf{r}_{\mathbf{k},\alpha}$, and the associated phase speed is $\omega_{\mathbf{k},B}=0$. The modes with $\alpha = \pm$ are the compressible gravity modes (also known as the Poincar\'e waves). The associated phase speed is $\omega_{\mathbf{k},\pm}= \pm \mbox{Ro}^{-1}\sqrt{|\mathbf{k}|^2 + 1}$,
and they are divergent.
Therefore,  the temporal evolution of these random Fourier coefficients satisfy
\begin{equation}\label{SI_SWE_Velocity_field_Fourier}
\begin{aligned}
\frac{\d\widehat{v}_{\mathbf{k},B}}{\d t} &= (-d_{\mathbf{k},B}+i\omega_{\mathbf{k},B}){\widehat{v}_{\mathbf{k},B}} + f_{\mathbf{k},B}(t) + \sigma_{\mathbf{k},B}\dot{W}_{\mathbf{k},B}(t),\\
\frac{ \d\widehat{v}_{\mathbf{k},\pm}}{\d t} &= (-d_{\mathbf{k},\pm}+i\omega_{\mathbf{k},\pm}){\widehat{v}_{\mathbf{k},\pm}} + f_{\mathbf{k},\pm}(t) + \sigma_{\mathbf{k},\pm}\dot{W}_{\mathbf{k},\pm}(t).
\end{aligned}
\end{equation}
The normalized eigenvector $\mathbf{r}_{\mathbf{k},B}$ of the GB modes are given by
\begin{equation}\label{SI_SWE_Modes0}
\mathbf{r}_{\mathbf{k},B} = \frac{1}{\sqrt{|\mathbf{k}|^2+1}}\left(
\begin{array}{c}
-ik_2 \\
ik_1 \\
1 \\
\end{array}
\right).
\end{equation}
The  normalized eigenvectors $\mathbf{r}_{\mathbf{k},\pm}$ of the corresponding to the gravity modes are given by
\begin{equation}\label{SI_SWE_Modespm}
\mathbf{r}_{\mathbf{k},\pm} = \frac{1}{|\mathbf{k}|\sqrt{(\delta + \delta^2)|\mathbf{k}|^2+2}}
\left(
\begin{array}{c}
ik_2\pm k_1\sqrt{\delta|\mathbf{k}|^2+1} \\
-ik_1\pm k_2\sqrt{\delta|\mathbf{k}|^2+1} \\
\delta |\mathbf{k}|^2 \\
\end{array}
\right).
\end{equation}
For the special case, $\mathbf{k}=\mathbf{0}$,
\begin{equation}\label{SI_SWE_Modespm0}
\mathbf{r}_{\mathbf{0},\pm} = \frac{1}{\sqrt{2}}
\left(
\begin{array}{c}
\pm i \\
1 \\
0 \\
\end{array}
\right).
\end{equation}

\subsection{Justification of Proposition \ref{theorem:random}}

	 Assume in total there are $L\gg 0$ tracers and the number of tracers for the randomly
	selected case is $0< L'<  L$.
	The time evolution of the covariance matrix for the two filters satisfy
	\begin{equation}\label{eq:evoR}
	\begin{aligned}
	\frac{\d\mathbf{R}_{\bf I}}{\d t} = &\boldsymbol{\Lambda}\mathbf{R}_{\bf I} + \mathbf{R}_{\bf I}\boldsymbol{\Lambda}^\ast + \boldsymbol{\Sigma}_\mathbf{U}\boldsymbol{\Sigma}_\mathbf{U}^\ast - \sigma_x^{-2}\mathbf{R}_{\bf I}\mathbf{P}_{\bf I}\mathbf{R}_{\bf I}\\
	\frac{\d\mathbf{R}_{\bf II}}{\d t}  = &\boldsymbol{\Lambda}\mathbf{R}_{\bf II} + \mathbf{R}_{\bf II}\boldsymbol{\Lambda}^\ast + \boldsymbol{\Sigma}_\mathbf{U}\boldsymbol{\Sigma}_\mathbf{U}^\ast - \sigma_x^{-2}\mathbf{R}_{\bf II}\mathbf{P}_{\bf II}\mathbf{R}_{\bf II}
	\end{aligned}
	\end{equation}
	where $\mathbf{P}_{\bf I}=\mathbf{A}_{\bf I}^\ast\mathbf{A}_{\bf I}$ with the entry being ${\bf P}_{{\bf I},mn}=\sum_{l=1}^{L}e^{i(m-n){\bf x}_l}\mathbf{r}^*_{{n}}\mathbf{r}_{{m}}$, while ${\bf P}_{\bf II}$ and ${\bf P}_{{\bf II},mn}$ have similar expressions.
	
	For the full filter, the mean-field theory gives an asymptotic limit of ${\bf R}_{\bf I}$, with the diagonal entry being given by \eqref{eq:diag_cov3}
	\begin{equation}\label{eq:asy_rI}
	r_{{\bf I},{\bf k}}=\frac{\sigma_{\bf k}^2}{d_{\bf k}+\sqrt{d_{\bf k}+L\sigma_x^{-2}\sigma_{\bf k}^2}}\sim\frac{\sigma_{\bf k}^2}{\sqrt{L}}
	\end{equation}
	For the approximate filter, ${\bf P}_{{\bf II},mn}=\sum_{l=1}^{L'}e^{i(m-n){\bf x}_l}\mathbf{r}^*_{{n}}\mathbf{r}_{{m}}$. Consider a small time interval
	{$[t, t+L\Delta t]$ } and use a time marching step $L\Delta t$ for a one-step discrete numerical approximation of integration of the evolution equation \eqref{eq:evoR}. Computing such an integration of the off-diagonal entries of ${\bf P_{II}}$ is equivalent to
	\begin{equation}\label{eq:PII}
	{\bf P}_{{\bf II}, mn}([t,t+L\Delta t])=\sum_{s=1}^{L}\sum_{l=1}^{L'}
	e^{t(m-n){\bf x}_l(t+s\Delta t)}{\bf r}_n^{\ast}{\bf r}_m
	=\sum_{l=1}^{LL'}e^{t(m-n){\bf x}_l(t)}{\bf r}_n^{\ast}{\bf r}_m+\mathcal{O}(\epsilon)
	\end{equation}
	where the assumption $\vert\mathbf{x}_l(t+L\Delta{t})-\mathbf{x}_l(t)\vert\leq \epsilon$ has been used. This, together with the discrete approximation of the numerical integration error accounts for the $O(\epsilon)$ term.
	Note that, in such an approximate filter, ${\bf x}_l$ at different time instants denotes different tracers. Therefore, making use of the facts that the tracers are uniformly distributed and $L$ is sufficiently large, the result in \eqref{eq:PII} implies the off-diagonal entries are of order
	${\bf P}_{\bf II}=\mathcal{O}(\epsilon)$.
	For the diagonal entries of $\bf R_{II}$, the mean-field theory again leads to
	\begin{equation}
	r_{{\bf II},{\bf k}}=\frac{\sigma_{\bf k}^2}{d_{\bf k}+\sqrt{d_{\bf k}+L'\sigma_x^{-2}\sigma_{\bf k}^2}}\sim\frac{\sigma_{\bf k}^2}{\sqrt{L'}}
	\end{equation}
	where $L'$ is assumed to be sufficiently large.
	
	Next, the time evolution of the posterior mean for the two filters are
	\begin{equation}\label{eq:evomu}
	\begin{aligned}
	\frac{\d\boldsymbol{\mu}_{\bf I}}{\d t} &= \left(\mathbf{F}_\mathbf{U} + \boldsymbol{\Lambda} \boldsymbol{\mu}_{\bf I}\right)  + \sigma_x^{-2}\mathbf{R}_{\bf I}\mathbf{A}_{\bf I}^\ast\left(\frac{\d \mathbf{X}}{\d t} - \mathbf{A}_{\bf I}\boldsymbol{\mu}_{\bf I} \right)\\
	\frac{\d\boldsymbol{\mu}_{\bf II}}{\d t} &= \left(\mathbf{F}_\mathbf{U} + \boldsymbol{\Lambda} \boldsymbol{\mu}_{\bf II}\right)  + \sigma_x^{-2}\mathbf{R}_{\bf II}\mathbf{A}_{\bf II}^{\ast}\left(\frac{\d \tilde{\mathbf{X}}}{\d t} - \mathbf{A}_{\bf II}\boldsymbol{\mu}_{\bf II} \right)
	\end{aligned}
	\end{equation}
	where
	\begin{equation}
	\mathbf{A}_{\bf I}=\begin{pmatrix}
	e^{i\bfk_1\bfx_1}{\bf r}_1 & e^{i\bfk_2\bfx_1}{\bf r}_2  & \cdots& e^{i\bfk_{|K|} \bfx_1}{\bf r}_{K}   \\
	e^{i\bfk_1\bfx_2}{\bf r}_1  & e^{i\bfk_2\bfx_2}{\bf r}_2  & \cdots& e^{i\bfk_{|K|} \bfx_2}{\bf r}_{K}  \\
	\vdots &\vdots &\ldots &\vdots\\
	e^{i\bfk_{1}\bfx_{L} }{\bf r}_1  & e^{i\bfk_{2}\bfx_{L}}{\bf r}_2   & \cdots& e^{i\bfk_{|K|}\bfx_L}{\bf r}_{K}  \\
	\end{pmatrix}\end{equation}
	and $\mathbf{A}_{\bf II}$ has a similar structure.

	In order to show the result \eqref{eq:random_obs_bounds}, namely $|\boldsymbol{\mu}_{\bf I}-\boldsymbol\mu_{\bf II}|\leq C\epsilon$, it is sufficient to compute the difference
	\begin{equation}\label{eq:musecondpart}
	\mathbf{R}_{\bf I}\mathbf{A}_{\bf I}^\ast\left(\frac{\d \mathbf{X}}{\d t} - \mathbf{A}_{\bf I}\boldsymbol{\mu}_{\bf I} \right) - \mathbf{R}_{\bf II}\mathbf{A}_{\bf II}^{\ast}\left(\frac{\d \tilde{\mathbf{X}}}{\d t} - \mathbf{A}_{\bf II}\boldsymbol{\mu}_{\bf II} \right)
	\end{equation}
	and see how if after certain manipulation it has an order of $\mathcal{O}(\epsilon)$.
	Apply a numerical approximation to integrate the time evolution of the posterior mean equation \eqref{eq:evomu} from $t=0$ to $t=L\Delta t$.
	For the full filter, the contribution from $\mathbf{A}_{\bf I}^\ast\frac{\d \mathbf{X}}{\d t}$ is
	\begin{equation}\label{eq:asy_muI_1}
	L\Delta t \mathbf{A}_{\bf I}^{\ast}\frac{\d \mathbf{X}}{\d t}=
	\begin{pmatrix}
	\sum_{l=1}^{L}e^{i\bfk_1\bfx_l}{\bf r}_1\frac{\d \mathbf{x}_l}{\d t}   \\
	\sum_{l=1}^{L}e^{i\bfk_2\bfx_l}{\bf r}_2\frac{\d \mathbf{x}_l}{\d t}  \\
	\vdots \\
	\sum_{l=1}^{L}e^{i\bfk_{|K|}\bfx_l}{\bf r}_{|K|}\frac{\d \mathbf{x}_l}{\d t}  \\
	\end{pmatrix}L\Delta t+\mathcal{O}(\epsilon)
	\end{equation}
	where the fact $\vert\mathbf{x}(t+L\Delta{t})-\mathbf{x}(t)\vert\leq \epsilon$ has been utilized. In addition, the mean-field theory indicates
	\begin{equation}\label{eq:asy_muI_2}
	L\Delta t{\bf A}_{\bf I}^{\ast}{\bf A}_{\bf I}\boldsymbol\mu_{\bf I}=L\boldsymbol\mu_{\bf I} L\Delta t+\mathcal{O}(\epsilon).
	\end{equation}
	Combining \eqref{eq:asy_rI}, \eqref{eq:asy_muI_1} and \eqref{eq:asy_muI_2} leads to
	\begin{equation}\label{eq:muI}
	\mathbf{R}_{\bf I}\mathbf{A}_{\bf I}^\ast\left(\frac{\d \mathbf{X}}{\d t} - \mathbf{A}_{\bf I}\boldsymbol{\mu}_{\bf I} \right)=\frac{1}{\sqrt{L}}L\Delta t
	\begin{pmatrix}
	\sum_{l=1}^{L}e^{i\bfk_1\bfx_l}{\bf r}_1\frac{\d \mathbf{x}_l}{\d t}-L\boldsymbol\mu_1   \\
	\sum_{l=1}^{L}e^{i\bfk_2\bfx_l}{\bf r}_2\frac{\d \mathbf{x}_l}{\d t}-L\boldsymbol\mu_2  \\
	\vdots \\
	\sum_{l=1}^{L}e^{i\bfk_{K}\bfx_l}{\bf r}_{K}\frac{\d \mathbf{x}_l}{\d t}-L\boldsymbol\mu_{K}  \\
	\end{pmatrix}+\mathcal{O}(\epsilon)
	\end{equation}
	
	For the approximate filter, similar results can be reached as the full filter case, except that $L$ should be replaced by $L'$,
	\begin{equation}
	L\Delta t \mathbf{A}_{\bf II}^\ast\frac{\d \tilde{\mathbf{X}}}{\d t}=
	\begin{pmatrix}
	\sum_{l=1}^{L'}e^{i\bfk_1\bfx_l}{\bf r}_1\frac{\d \mathbf{x}_l}{\d t}   \\
	\sum_{l=1}^{L'}e^{i\bfk_2\bfx_l}{\bf r}_2\frac{\d \mathbf{x}_l}{\d t}  \\
	\vdots \\
	\sum_{l=1}^{L'}e^{i\bfk_{K}\bfx_l}{\bf r}_{K}\frac{\d \mathbf{x}_l}{\d t}  \\
	\end{pmatrix}L\Delta t+\mathcal{O}(\epsilon).
	\end{equation}
	Utilizing the fact that, at different time instants, the tracers
	are randomly selected and are different with each other, it can be shown that
	$L\Delta t{\bf A}_{\bf II}^{\ast}{\bf A}_{\bf II}\boldsymbol\mu_{\bf II}=L'\boldsymbol\mu_{\bf II} L\Delta t$.
	Therefore, after integrating $L$ steps, the following result holds
	\begin{equation}\label{eq:muII}
	\mathbf{R}_{\bf II}\mathbf{A}_{\bf II}^\ast\left(\frac{\d \tilde{\mathbf{X}}}{\d t} - \mathbf{A}_{\bf II}\boldsymbol{\mu}_{\bf II} \right)=\frac{1}{\sqrt{L'}}L'\Delta t
	\begin{pmatrix}
	\sum_{l=1}^{L}e^{i\bfk_1\bfx_l}{\bf r}_1\frac{\d \mathbf{x}_l}{\d t}-L\mu_1   \\
	\sum_{l=1}^{L}e^{i\bfk_2\bfx_l}{\bf r}_2\frac{\d \mathbf{x}_l}{\d t}-L\mu_2  \\
	\vdots \\
	\sum_{l=1}^{L}e^{i\bfk_{K}\bfx_l}{\bf r}_{K}\frac{\d \mathbf{x}_l}{\d t}-L\mu_{K}  \\
	\end{pmatrix}+\mathcal{O}(\epsilon)
	\end{equation}
	By comparing \eqref{eq:muI} and \eqref{eq:muII}, it is seen that, after multiplying a perfactor $\frac{\sqrt{L}}{\sqrt{L'}}$ to \eqref{eq:muII},
	the result in \eqref{eq:musecondpart} will have an order of $\mathcal{O}(\epsilon)$.

	\bibliographystyle{plain}
\bibliography{references}
\end{document}